\documentclass[10pt,a4paper]{amsart}
\usepackage{amsthm,amssymb,amsmath,mathabx}
\usepackage{amsfonts,enumerate}
\usepackage{mathrsfs}
\usepackage{graphicx}
\usepackage{xypic}
\usepackage[body={6.0in,8.2in},left=1.00in,right=1.00in]{geometry}

\theoremstyle{plain}
\newtheorem{thm}{Theorem}[section]
\newtheorem*{thm*}{Theorem}
\newtheorem{prop}[thm]{Proposition}
\newtheorem*{prop*}{Proposition}
\newtheorem{lemma}[thm]{Lemma}
\newtheorem*{lemma*}{Lemma}
\newtheorem{coro}[thm]{Corollary}

\theoremstyle{definition}

\newtheorem*{dfn*}{Definition}

\theoremstyle{remark}
\newtheorem{rem}[thm]{Remark}
\newtheorem*{rem*}{Remark}
\newtheorem*{rems*}{Remarks}
\newtheorem*{fact*}{Fact}

\newtheorem*{exc*}{Exercise}

\newcommand{\A}{\mathbb A}
\newcommand{\C}{\mathbb C}
\newcommand{\R}{\mathbb R}

\newcommand{\Z}{\mathbb Z}
\newcommand{\F}{\mathbb F}
\newcommand{\ba}{\mathbf a}
\newcommand{\bbc}{\mathbf b}
\newcommand{\Q}{\mathbb Q}

\newcommand{\m}{\mathfrak m}

\newcommand{\Or}{\mathscr O}
\newcommand{\Gal}{\operatorname{Gal}}

\newcommand{\Art}{\operatorname{Art}}
\newcommand{\BC}{\operatorname{BC}}
\newcommand{\Mod}{\operatorname{mod}}
\newcommand{\tr}{\operatorname{Tr}}
\newcommand{\ssp}{\operatorname{ss}}

\newcommand{\HT}{\operatorname{HT}}
\newcommand{\Hom}{\operatorname{Hom}}
\newcommand{\Aut}{\operatorname{Aut}}

\newcommand{\opp}{\operatorname{opp}}

\newcommand{\End}{\operatorname{End}}
\newcommand{\Res}{\operatorname{Res}}
\newcommand{\GL}{\operatorname{GL}}

\newcommand{\T}{\mathbb{T}}

\newcommand{\vabs}{|\hspace{1mm}|}
\newcommand{\ad}{\operatorname{ad}}
\newcommand{\red}{\operatorname{red}}
\newcommand{\Frob}{\operatorname{Frob}}
\newcommand{\depth}{\operatorname{depth}}

\newcommand{\ab}{\operatorname{ab}}

\newcommand{\gr}{\operatorname{gr}}
\newcommand{\loc}{\operatorname{loc}}

\newcommand{\diag}{\operatorname{diag}}
\newcommand{\ind}{\operatorname{ind}}
\newcommand{\univ}{\operatorname{univ}}
\newcommand{\lin}{\operatorname{lin}}

\newcommand{\Iw}{\operatorname{Iw}}

\newcommand{\Sp}{\operatorname{Sp}}
\newcommand{\St}{\operatorname{St}}
\newcommand{\Ann}{\operatorname{Ann}}

\newcommand{\G}{\mathscr{G}}
\newcommand{\CS}{\mathscr{S}}

\newcommand{\dR}{\operatorname{dR}}

\numberwithin{equation}{subsection}

\begin{document} 

\title{Modularity lifting theorems for Galois representations of unitary type}

\author{Lucio Guerberoff}
\email{lguerb@math.jussieu.fr}
\address{Institut de Math\'ematiques de Jussieu, Universit\'e Paris 7 Denis Diderot, 175 rue du Chevaleret, 75013, Paris, France - Departamento de Matem\'atica, Universidad de Buenos Aires, Pabell\'on I, Ciudad Universitaria, 1428, Buenos Aires, Argentina}


\thanks{The author was partially supported by a CONICET fellowship and a fellowship from the Ministerio de Educaci\'on de Argentina \& Ambassade de France en Argentine}

\begin{abstract}
 We prove modularity lifting theorems for $\ell$-adic Galois representations of any dimension satisfying a unitary type condition and a Fontaine-Laffaille condition at $\ell$. This extends the results of Clozel, Harris and Taylor, and the subsequent work by Taylor. The proof uses the Taylor-Wiles method, as improved by Diamond, Fujiwara, Kisin and Taylor, applied to Hecke algebras of unitary groups, and results of Labesse on stable base change and descent from unitary groups to $\GL_n$.
\end{abstract}

\maketitle

\setcounter{section}{-1}

\section*{Introduction} The goal of this paper is to prove modularity lifting theorems for Galois representations of any dimension satisfying certain conditions. We largely follow the articles \cite{cht} and \cite{taylorII}, where an extra local condition appears. In this work we remove that condition, which can be done thanks to the latest developments of the trace formula.  More precisely, let $F$ be a totally imaginary quadratic extension of a totally real field $F^+$. Let $\Pi$ be a cuspidal automorphic representation of $\GL_n(\A_F)$ satisfying the following conditions.
\begin{itemize}
\item There exists a continuous character $\chi:\A_{F^{+}}^{\times}/(F^{+})^{\times}\to\C^{\times}$ such that $\chi_{v}(-1)$ is independent of $v|\infty$ and
\[ \Pi^\vee\cong\Pi^{c}\otimes(\chi\circ\mathbf{N}_{F/F^{+}}\circ\det).\]
\item $\Pi$ is cohomological.
\end{itemize}
Here, $c$ is the non-trivial Galois automorphism of $F/F^+$, and {\em cohomological} means that $\Pi_\infty$ has the same infinitesimal character as an algebraic, finite dimensional, irreducible representation of $(\Res_{F/\Q}\GL_n)(\C)$. Let $\ell$ be a prime number, and $\iota:\overline\Q_\ell\overset{\sim}{\longrightarrow}\C$ an isomorphism. Then there is a continuous semisimple Galois representation
\[ r_{\ell,\iota}(\Pi):\Gal(\overline F/F)\to\GL_n(\overline\Q_\ell) \]
which satisfies certain expected conditions. In particular, for places $v$ of $F$ not dividing $\ell$, the restriction $r_{\ell,\iota}(\Pi)|_{\Gal(\overline F_v/F_v)}$ to a decomposition group at $v$ should be isomorphic, as a Weil-Deligne representation, to the representation corresponding to $\Pi_v$ under a suitably normalized local Langlands correspondence. The construction of the 
Galois representation $r_{\ell,\iota}(\Pi)$ under these hypotheses is due to Clozel, Harris and Labesse (\cite{chl,chl1}), Chenevier and Harris (\cite{chenharr}), and Shin
(\cite{shin}), although they only match the Weil parts and not the whole Weil-Deligne representation. In the case that $\Pi$ satisfies the additional hypothesis that $\Pi_{v}$ is a square integrable representation for some finite place $v$, Taylor and Yoshida have shown in \cite{ty} that the corresponding Weil-Deligne representations are indeed the same, as expected. Without the square integrable hypothesis, this is proved by Shin in \cite{shin} in the case where $n$ is odd, or when $n$ is even and the archimedean weight of $\Pi$ is 'slightly regular', a mild condition we will not recall here. We will not need this stronger result for the purposes of our paper.

We use the instances of stable base change and descent from $\GL_n$ to unitary groups, proved by Labesse (\cite{lab}) to attach Galois representations to automorphic representations of totally definite unitary groups. In this setting, we prove an $R^{\red}=T$ theorem, following the development of the Taylor-Wiles method used in \cite{taylorII}. Finally, using the results of Labesse again, we prove our modularity lifting theorem for $\GL_n$. We describe with more detail the contents of this paper.

Section 1 contains some basic preliminaries. We include some generalities about smooth representations of $\GL_n$ of a $p$-adic field, over $\overline\Q_\ell$ or $\overline\F_\ell$, which will be used later in the proof of the main theorem. We note that many of the results of this section are also proved in \cite{cht}, although in a slightly different way. We stress the use of the Bernstein formalism in our proofs; some of them are based on an earlier draft \cite{4thdraft} of \cite{cht}.

In Section 2, we develop the theory of ($\ell$-adic) automorphic forms on totally definite unitary groups, and apply the results of Labesse and the construction mentioned above to attach Galois representations to automorphic representations of unitary groups. 

In Section 3, we study the Hecke algebras of unitary groups and put everything together to prove the main result of the paper. More precisely, if $\T$ denotes the (localized) Hecke algebra and $R$ is a certain universal deformation ring of a mod $\ell$ Galois representation attached to $\T$, we prove that $R^{\red}=\T$. In Section 4, we go back to $\GL_n$ and use this result to prove the desired modularity lifting theorems. The most general theorem we prove for imaginary CM fields is the following. For the terminology used in the different hypotheses, we refer the reader to the main text. 

\begin{thm*} Let $F^+$ be a totally real field, and $F$ a totally imaginary quadratic extension of $F^+$. Let $n\geq 1$ be an integer and $\ell>n$ be a prime number, unramified in $F$. Let
 \[ r:\Gal(\overline F/F)\longrightarrow\GL_n(\overline\Q_\ell) \]
be a continuous irreducible representation with the following properties. Let $\overline r$ denote the semisimplification of the reduction of $r$.

\begin{enumerate}[(i)]
 \item $r^c\cong r^\vee(1-n)$.
 \item $r$ is unramified at all but finitely many primes.
 \item For every place $v|\ell$ of $F$, $r|_{\Gamma_v}$ is crystalline.
  \item There is an element $\ba\in(\Z^{n,+})^{\Hom(F,\overline\Q_\ell)}$ such that

     \begin{itemize} 
         \item for all $\tau\in\Hom(F^+,\overline\Q_\ell)$, we have either
		\[ \ell-1-n\geq a_{\tau,1}\geq\dots\geq a_{\tau,n}\geq 0 \]
		or
		\[ \ell-1-n\geq a_{\tau c,1}\geq\dots\geq a_{\tau c,n}\geq 0;\]
	\item for all $\tau\in\Hom(F,\overline\Q_\ell)$ and every $i=1,\dots,n$,
    		 \[ a_{\tau c,i}=-a_{\tau,n+1-i}. \]
	\item for all $\tau\in\Hom(F,\overline\Q_\ell)$ giving rise to a prime $w|\ell$,
	\[ \HT_{\tau}(r|_{\Gamma_w})=\{j-n-a_{\tau,j}\}_{j=1}^n.\]
	In particular, $r$ is Hodge-Tate regular.
	\end{itemize}
 \item $\overline F^{\ker(\ad\overline r)}$ does not contain $F(\zeta_\ell)$.
 \item The group $\overline r(\Gal(\overline F/F(\zeta_\ell)))$ is big.
 \item The representation $\overline r$ is irreducible and there is a conjugate self-dual, cohomological, cuspidal automorphic representation $\Pi$ of $\GL_{n}(\A_{F})$, of weight $\ba$ and unramified above $\ell$, and an isomorphism $\iota:\overline\Q_\ell\overset{\sim}{\longrightarrow}\C$, such that $\overline r\cong\overline r_{\ell,\iota}(\Pi)$.
\end{enumerate}
Then $r$ is automorphic of weight $\ba$ and level prime to $\ell$.

\end{thm*}

We make some remarks about the conditions in the theorem. Condition (i) says that $r$ is conjugate self-dual, and this is essential for the numerology behind the Taylor-Wiles method. Conditions (ii) and (iii) say that the Galois representation is geometric in the sense of Fontaine-Mazur, although it says a little more. It is expected that one can relax condition (iii) to the requirement that $r$ is de Rham at places dividing $\ell$. The stronger crystalline form, the hypothesis on the Hodge-Tate weights made in (iv) and the requirement that $\ell>n$ is unramified in $F$ are needed to apply the theory of Fontaine and Laffaille to calculate the local deformation rings. The condition that $\ell>n$ is also used to treat non-minimal deformations. Condition (v) allows us to choose auxiliary primes to augment the level and ensure that certain level structures are sufficiently small. The bigness condition in (vi) is to make the Tchebotarev argument in the Taylor-Wiles method work. Hypothesis (vii) is, as usual, essential to the method. An analogous theorem can be proved over totally real fields.

{\bf Acknowledgements.}I take this opportunity to express my gratitude to my thesis adviser Michael Harris for his many explanations and suggestions, and specially for his advices. I would also like to thank Brian Conrad, Jean-Francois Dat, Alberto M\'inguez and Nicol\'as Ojeda B\"ar for useful conversations. Finally, I would like to thank Roberto Miatello for his constant support and encouragement. I also thank the referee for his corrections and helpful comments.

\section{Some notation and definitions}
As a general principle, whenever $F$ is a field and $\overline F$ is a chosen separable closure, we write $\Gamma_F=\Gal(\overline F/F)$. We also write $\Gamma_F$ when the choice of $\overline F$ is implicit. If $F$ is a number field and $v$ is a place of $F$, we usually write $\Gamma_v\subset\Gamma_F$ for a decomposition group at $v$. If $v$ is finite, we denote by $q_v$ the order of the residue field of $v$.

\subsection{Irreducible algebraic representations of $\GL_n$}\label{irrgln} Let $\Z^{n,+}$ denote the set of $n$-tuples of integers $a=(a_1,\dots,a_n)$ such that 
\[ a_1\geq\dots\geq a_n.\]
Given $a\in\Z^{n,+}$, there is a unique irreducible, finite dimensional, algebraic representation $\xi_a:\GL_n\to\GL(W_a)$ over $\Q$ with highest weight
given by
\[ \diag(t_1,\dots,t_n)\mapsto\prod_{i=1}^nt_i^{a_i}.\]
Let $E$ be any field of characteristic zero. Tensoring with $E$, we obtain an irreducible algebraic representation $W_{a,E}$ of $\GL_n$ over $E$, and every such representation arises in
this way. Suppose that $E/\Q$ is a finite extension. Then the irreducible, finite dimensional, algebraic representations of $(\Res_{E/\Q}\GL_{n/E})(\C)$ are parametrized by elements $\ba\in(\Z^{n,+})^{\Hom(E,\C)}$. We denote them by $(\xi_\ba,W_\ba)$.

\subsection{Local Langlands correspondence}
Let $p$ be a rational prime and let $F$ be a finite extension of $\Q_p$. Fix an algebraic closure $\overline F$ of $F$. Fix also a positive integer $n$, a prime number $\ell\neq p$ and an algebraic closure $\overline\Q_\ell$ of $\Q_\ell$. Let $\Art_{F}:F^\times\to \Gamma_F^{\ab}$ be the local reciprocity map, normalized to take uniformizers to geometric Frobenius elements. If $\pi$ is an irreducible smooth representation of $\GL_n(F)$ over $\overline\Q_\ell$, we will write $r_\ell(\pi)$ for the $\ell$-adic Galois representation associated to the Weil-Deligne representation
\[ \mathscr{L}(\pi\otimes\vabs^{(1-n)/2}),\]
where $\mathscr{L}$ denotes the local Langlands correspondence, normalized to coincide with the correspondence induced by $\Art_{F}$ in the case $n=1$. Note that $r_\ell(\pi)$ does not always exist. The eigenvalues of $\mathscr{L}(\pi\otimes\vabs^{(n-1)/2})(\phi_F)$ must be $\ell$-adic units for some lift $\phi_F$ of the geometric Frobenius (see \cite{tate}). Whenever we make a statement about $r_\ell(\pi)$, we will suppose that this is the case. Note that our conventions differ from those of \cite{cht} and \cite{taylorII}, where $r_{\ell}(\pi)$ is defined to be the Galois representation associated to $\mathscr{L}(\pi^{\vee}\otimes\vabs^{(1-n)/2})$.

\subsection{Hodge-Tate weights}
Fix a finite extension $L/\Q_\ell$ and an algebraic closure $\overline L$ of $L$. Fix an algebraic closure $\overline\Q_\ell$ of $\Q_{\ell}$ and an algebraic extension $K$ of $\Q_\ell$ contained in $\overline\Q_\ell$ such that $K$ contains every $\Q_\ell$-embedding $L\hookrightarrow\overline\Q_\ell$. Suppose that $V$ is a finite dimensional $K$-vector space equipped with a continuous linear action of $\Gamma_L$. Let $B_{\dR}$ be the ring of $p$-adic periods, as in \cite{Ast}. Then $(B_{\dR}\otimes_{\Q_\ell}V)^{\Gamma_L}$ is an $L\otimes_{\Q_\ell}K$-module. We say that $V$ is de Rham if this module is free of rank equal to $\dim_{K}V$. Since $L\otimes_{\Q_\ell}K\simeq(K)^{\Hom_{\Q_\ell}(L,K)}$, if $V$ is a $K$-representation of $\Gamma_L$, we have that
\[
\begin{array}{rl}
 (B_{\dR}\otimes_{\Q_\ell}V)^{\Gamma_L} & \simeq\underset{\tau\in\Hom_{\Q_\ell}(L,K)}\prod (B_{\dR}\otimes_{\Q_\ell}V)^{\Gamma_L}\otimes_{L\otimes_{\Q_\ell}K,\tau\otimes 1}K \\
 & \simeq\underset{\tau\in\Hom_{\Q_\ell}(L,K)}\prod(B_{\dR}\underset{L,\tau}\otimes V)^{\Gamma_L}.\end{array}\]
It follows that $V$ is de Rham if and only if
\[ \dim_{K}(B_{\dR}\otimes_{L,\tau}V)^{\Gamma_L}=\dim_{K}V \]
for every $\tau\in\Hom_{\Q_\ell}(L,K)$. We use the convention of Hodge-Tate weights in which the cyclotomic character has $1$ as its unique Hodge-Tate weight. Thus, for $V$ de Rham, we let $\HT_\tau(V)$ be the multiset consisting of the elements $q\in\Z$ such that $\gr^{-q}(B_{\dR}\otimes_{L,\tau}V)^{\Gamma_L}\neq 0$, with multiplicity equal to
\[ \dim_{K}\gr^{-q}(B_{\dR}\otimes_{L,\tau}V)^{\Gamma_L}.\]
Thus, $\HT_\tau(V)$ is a multiset of $\dim_{K}V$ elements. We say that $V$ is {\em Hodge-Tate regular} if for every $\tau\in\Hom_{\Q_\ell}(L,K)$, the multiplicity of each Hodge-Tate weight with respect to $\tau$ is $1$. We make analogous definitions for crystalline representations over $K$.

\subsection{Galois representations of unitary type}\label{unittype} Let $F$ be any number field. If $\ell$ is a prime number, $\iota:\overline\Q_\ell\overset{\sim}{\longrightarrow}\C$ is an isomorphism and $\psi:\A_F^\times/F^\times\to\C^\times$ is an algebraic character, we denote by $r_{\ell,\iota}(\psi)$ the Galois character associated to it by Lemma 4.1.3 of \cite{cht}.

Let $F^+$ be a totally real number field, and $F/F^+$ a totally imaginary
quadratic extension. Denote by $c\in\Gal(F/F^+)$ the non-trivial automorphism. Let $\Pi$ be an irreducible admissible representation of $\GL_n(\A_F)$. We say that $\Pi$ is {\em essentially conjugate self dual} if there exists a continuous character $\chi:\A_{F^{+}}^{\times}/(F^{+})^{\times}\to\C^{\times}$ with $\chi_{v}(-1)$ independent of $v|\infty$ such that
\[ \Pi^\vee\cong\Pi^{c}\otimes(\chi\circ\mathbf{N}_{F/F^{+}}\circ\det).\]
If we can take $\chi=1$, that is, if $\Pi^{\vee}\cong\Pi^{c}$, we say that $\Pi$ is {\em conjugate self dual}.

Let $\Pi$ be an automorphic representation of $\GL_n(\A_F)$. We say that $\Pi$ is {\em cohomological} if there exists an irreducible, algebraic, finite-dimensional representation $W$ of $\Res_{F/\Q}\GL_n$, such that the infinitesimal character of $\Pi_\infty$ is the same as
that of $W$. Let $\ba\in(\Z^{n,+})^{\Hom(F,\C)}$, and let $(\xi_\ba,W_\ba)$ the irreducible, finite dimensional, algebraic representation of $(\Res_{F/\Q}\GL_n)(\C)$ with highest weight $\ba$. We say that $\Pi$ has {\em weight} $\ba$ if it has the same infinitesimal character as $(\xi_\ba^\vee,W_\ba^\vee)$. 

The next theorem (in the conjugate self dual case) is due to Clozel, Harris and Labesse (\cite{chl,chl1}), with some improvements by Chenevier and Harris (\cite{chenharr}), except that they only provide compatibility of the local and global Langlands correspondences for the unramified places. Shin (\cite{shin}), using a very slightly different method, obtained the identification at the remaining places. The slightly more general version stated here for an essentially conjugate self dual representation is proved in Theorem 1.2 of \cite{cy2}. Let $\overline F$ be an algebraic closure of $F$ and let $\Gamma_F=\Gal(\overline F/F)$. For $m\in\Z$ and $r:\Gamma_F\to\GL_n(\overline\Q_\ell)$ a continuous representation, we denote by $r(m)$ the $m$-th Tate twist of $r$, and by $r^{\ssp}$ the semisimplification of $r$. Fix a prime number $\ell$, an algebraic closure $\overline\Q_\ell$ of $\Q_\ell$, and an isomorpshim $\iota:\overline\Q_\ell\overset{\sim}{\longrightarrow}\C$.

\begin{thm} Let $\Pi$ be an essentially conjugate self dual, cohomological, cuspidal automorphic representation of $\GL_n(\A_F)$. More precisely, suppose that $\Pi^{\vee}\cong\Pi^{c}\otimes(\chi\circ\mathbf{N}_{F/F^{+}}\circ\det)$ for some continuous character $\chi:\A_{F^{+}}^{\times}/(F^{+})^{\times}\to\C^{\times}$ with $\chi_{v}(-1)$ independent of $v|\infty$. Then there exists a continuous semisimple representation
 \[ r_\ell(\Pi)=r_{\ell,\iota}(\Pi):\Gamma_F\to\GL_n(\overline\Q_\ell) \]
with the following properties.
\begin{enumerate}[(i)]
 \item For every finite place $w\nmid\ell$,
 \[ (r_\ell(\Pi)|_{\Gamma_w})^{\ssp}\simeq\left(r_\ell(\iota^{-1}\Pi_w)\right)^{\ssp}.\]
 \item $r_\ell(\Pi)^c\cong r_\ell(\Pi)^\vee(1-n)\otimes r_{\ell}(\chi^{-1})|_{\Gamma_{F}}$.
 \item If $w\nmid\ell$ is a finite place such that $\Pi_w$ is unramified, then $r_\ell(\Pi)$ is unramified at $w$.
 \item For every $w|\ell$, $r_\ell(\Pi)$ is de Rham at $w$. Moreover, if $\Pi_w$ is unramified, then $r_\ell(\Pi)$ is crystalline at $w$.
 \item Suppose that $\Pi$ has weight $\ba$. Then for each $w|\ell$ and each embedding $\tau:F\hookrightarrow\overline\Q_\ell$ giving rise to $w$, the Hodge-Tate weights of $r_\ell(\Pi)|_{\Gamma_w}$ with respect to $\tau$ are given by
 \[ \HT_\tau(r_\ell(\Pi)|_{\Gamma_w})=\{j-n-a_{\iota\tau,j}\}_{j=1,\dots,n},\]
 and in particular, $r_\ell(\Pi)|_{\Gamma_w}$ is Hodge-Tate regular.
\end{enumerate}

\end{thm}

The representation $r_{\ell,\iota}(\Pi)$ can be taken to be valued in the ring of integers of a finite extension of $\Q_\ell$. Thus, we can reduce it modulo its maximal ideal and semisimplify to obtain a well defined continuous semisimple representation
\[ \overline r_{\ell,\iota}(\Pi):\Gamma_F\longrightarrow\GL_n(\overline\F_\ell).\]

Let $\ba$ be an element of $(\Z^{n,+})^{\Hom(F,\overline\Q_\ell)}$. Let
\[ r:\Gamma_F\longrightarrow\GL_n(\overline\Q_\ell)
\]
be a continuous semisimple representation. We say that $r$ is {\em automorphic of weight} $\ba$ if there is an isomorphism $\iota:\overline\Q_\ell\overset\sim\longrightarrow\C$ and an essentially conjugate self dual, cohomological, cuspidal automorphic representation $\Pi$ of $\GL_n(\A_F)$ of weight $\iota_*\ba$ such that $r\cong r_{\ell,\iota}(\Pi)$. We say that $r$ is automorphic of weight $\ba$ and {\em level prime to} $\ell$ if moreover there exists such a pair $(\iota,\Pi)$ with $\Pi_\ell$ unramified. Here $\iota_*\ba\in(\Z^{n,+})^{\Hom(F,\C)}$ is defined as $(\iota_*\ba)_\tau=a_{\iota^{-1}\tau}$.

There is an analogous construction for a totally real field $F^+$. The definition of cohomological is the same, namely, that the infinitesimal character is the same as that of some irreducible algebraic finite dimensional representation of $(\Res_{F^+/\Q}\GL_n)(\C)$.

\begin{thm} Let $\Pi$ be a cuspidal automorphic representation of $\GL_n(\A_{F^+})$, cohomological of weight $\ba$, and suppose that
\[ \Pi^\vee\cong\Pi\otimes(\chi\circ\det),\]
where $\chi:\A_{F^+}^\times/(F^+)^\times\to\C^\times$ is a continuous character such that $\chi_v(-1)$ is independent of $v|\infty$. Let $\iota:\overline\Q_\ell\overset{\sim}{\longrightarrow}\C$. Then there is a continuous semisimple representation
\[ r_{\ell}(\Pi)=r_{\ell,\iota}(\Pi):\Gamma_{F^+}\to\GL_n(\overline\Q_\ell) \]
with the following properties.

\begin{enumerate}[(i)]
 \item For every finite place $v\nmid\ell$,
 \[ (r_\ell(\Pi)|_{\Gamma_v})^{\ssp}\simeq\left(r_\ell(\iota^{-1}\Pi_v)\right)^{\ssp}.\]
 \item $r_\ell(\Pi)\cong r_\ell(\Pi)^\vee(1-n)\otimes r_\ell(\chi^{-1})$.
 \item If $v\nmid\ell$ is a finite place such that $\Pi_v$ is unramified, then $r_\ell(\Pi)$ is unramified at $v$.
 \item For every $v|\ell$, $r_\ell(\Pi)$ is de Rham at $v$. Moreover, if $\Pi_v$ is unramified, then $r_\ell(\Pi)$ is crystalline at $v$.
 \item For each $v|\ell$ and each embedding $\tau:F^+\hookrightarrow\overline\Q_\ell$ giving rise to $v$, the Hodge-Tate weights of $r_\ell(\Pi)|_{\Gamma_v}$ with respect to $\tau$ are given by
 \[ \HT_\tau(r_\ell(\Pi)|_{\Gamma_v})=\{j-n-a_{\iota\tau,j}\}_{j=1,\dots,n},\]
 and in particular, $r_\ell(\Pi)|_{\Gamma_v}$ is Hodge-Tate regular.
\end{enumerate}
Moreover, if $\psi:\A_{F^+}^\times/(F^+)^\times\to\C^\times$ is an algebraic character, then
\[ r_\ell(\Pi\otimes(\psi\circ\det))=r_\ell(\Pi)\otimes r_\ell(\psi).\]

\begin{proof} This can be deduced from the last theorem in exactly the same way as Proposition 4.3.1 of \cite{cht} is deduced from Proposition 4.2.1 of {\em loc. cit.}
\end{proof}
\end{thm}

We analogously define what it means for a Galois representation of a totally real field to be automorphic of some weight $\ba$.

\section{Admissible representations of $\GL_n$ of a $p$-adic field over $\overline\Q_\ell$ and $\overline\F_\ell$}\label{localtheory}
Let $p$ be a rational prime and let $F$ be a finite extension of $\Q_p$, with ring of integers $\Or_F$, maximal ideal $\lambda_F$ and residue field $k_F=\Or_F/\lambda_F$. Let $q=\#k_F$. Let $\overline\omega$ be a generator of $\lambda_F$. We will fix an algebraic closure $\overline F$ of $F$, and write $\Gamma_F=\Gal(\overline F/F)$. Corresponding to it, we have an algebraic closure $\overline{k_F}$ of $k_F$, and we will let $\Frob_F$ be the geometric Frobenius in $\Gal(\overline{k_F}/k_F)$ and $I_F$ be the inertia subgroup of $\Gamma_F$. Usually we will also write $\Frob_F$ for a lift to $\Gamma_F$. Fix also a positive integer $n$, a prime number $\ell\neq p$, an algebraic closure $\overline\Q_\ell$ of $\Q_\ell$ and an algebraic closure $\overline\F_\ell$ of $\F_\ell$. We will let $R$ be either $\overline\Q_\ell$ or $\overline\F_\ell$. Denote by $\vabs:F^{\times}\to q^{\Z}\subset\Z[\frac{1}{q}]$ the absolute value normalized such that $|\overline\omega|=q^{-1}$. We denote by the same symbol the composition of $\vabs$ and the natural map $\Z[\frac{1}{q}]\to R$, which exists because $q$ is invertible in $R$. For the general theory of smooth representations over $R$, we refer the reader to \cite{vigirr}. Throughout this section, representation will always mean smooth representation.

For a locally compact, totally disconnected group $G$, a compact open subgroup $K\subset G$ and an element $g\in G$, we denote by $[KgK]$ the operator in the Hecke algebra of $G$ relative to $K$ corresponding to the ($R$-valued) characteristic function of the double coset $KgK$.

Given a tuple $\mathbf{t}=(t^{(1)},\dots,t^{(n)})$ of elements in any ring $A$, we denote by $P_{q,\mathbf{t}}\in A[X]$ the polynomial
        \[ P_{q,\mathbf{t}}=X^n+\sum_{j=1}^n(-1)^jq^{j(j-1)/2}t^{(j)}X^{n-j}. \]        

We use freely the terms Borel, parabolic, Levi, and so on, to refer to the $F$-valued points of the corresponding algebraic subgroups of $\GL_n$. Write $B$ for the Borel subgroup of $\GL_n(F)$ consisting of upper triangular matrices, and $B_0=B\cap\GL_n(\Or_F)$. Let $T\simeq (F^\times)^n$ be the standard maximal torus of $\GL_n(F)$. Let $N$ be the group of upper triangular matrices whose diagonal elements are all $1$. Then $B=TN$ (semi-direct product). Let $r:\GL_n(\Or_F)\to\GL_n(k_F)$ denote the reduction map. We introduce the following subgroups of $\GL_n(\Or_F)$:

\begin{itemize}
\item $U_0=\{g\in\GL_n(\Or_F):\hspace{2mm}r(g)=\left(\begin{array}{cc}*_{n-1,n-1}&*_{n-1,1}\\0_{1,n-1}&*\end{array}\right)\}$;
\item $U_1=\{g\in\GL_n(\Or_F):\hspace{2mm}r(g)=\left(\begin{array}{cc}*_{n-1,n-1}&*_{n-1,1}\\0_{1,n-1}&1\end{array}\right)\}$;
\item $\Iw=\{g\in\GL_n(\Or_F):\hspace{2mm}r(g)\text{ is upper triangular}\}$;
\item $\Iw_1=\{g\in\Iw:\hspace{2mm}r(g)_{ii}=1\hspace{1mm}\forall i=1,\dots,n\}$.
\end{itemize}

Thus, $U_1$ is a normal subgroup of $U_0$ and we have a natural identification
\[ U_0/U_1\simeq k_F^\times,\]
and similarly, $\Iw_1$ is a normal subgroup of $\Iw$ and we have a natural identification
\[ \Iw/\Iw_1\simeq(k_F^\times)^n.\]

We denote by $\mathscr{H}$ the $R$-valued Hecke algebra of $\GL_n(F)$ with respect to $\GL_n(\Or_F)$. We do not include $R$ in the notation. For every smooth representation $\pi$ of $\GL_n(F)$, $\pi^{\GL_n(\Or_F)}$ is naturally a left module over $\mathscr{H}$. For $j=1,\dots,n$, we will let $T_F^{(j)}\in\mathscr{H}$ denote the Hecke operator
\[ \left[\GL_n(\Or_F)\left(\begin{array}{cc}\overline\omega1_j&0\\0&1_{n-j}\end{array}\right)\GL_n(\Or_F)\right].\]

Let $\pi$ be a representation of $\GL_n(F)$ over $\overline\Q_\ell$. We say that $\pi$ is essentially square-integrable if, under an isomorphism $\overline\Q_\ell\cong\C$, the corresponding complex representation is essentially square integrable in the usual sense. It is a non trivial fact that the notion of essentially square integrable complex representation is invariant under an automorphism of $\C$, which makes our definition independent of the chosen isomorphism $\overline\Q_\ell\cong\C$. This can be shown using the Bernstein-Zelevinsky classification of essentially square integrable representations in terms of quotients of parabolic inductions from supercuspidals (see below).

Let $n=n_1+\dots+n_r$ be a partition of $n$ and $P\supset B$ the corresponding parabolic subgroup of $\GL_n(F)$. The modular character $\delta_P:P\to\Q^\times$ takes values in $q^{\Z}\subset R^{\times}$. Choosing once and for all a square root of $q$ in $R$, we can consider the square root character $\delta_{P}^{1/2}:P\to R^{\times}$. For each $i=1,\dots r$, let $\pi_i$ be a representation of $\GL_{n_i}(F)$. We denote by $\pi_1\times\dots\times\pi_r$ the normalized induction from $P$ to $\GL_n(F)$ of the representation $\pi_1\otimes\dots\otimes\pi_r$. Whenever we write $\vabs$ we will mean $\vabs\circ\det$. For any $R$-valued character $\beta$ of $F^\times$ and any positive integer $m$, we denote by $\beta[m]$ the one dimensional representation $\beta\circ\det$ of $\GL_m(F)$.

Suppose that $R=\overline\Q_\ell$. Let $n=rk$ and $\sigma$ be an irreducible supercuspidal representation of $\GL_r(F)$. By a theorem of Bernstein (\cite[9.3]{zelevII}), 
\[ \left(\sigma\otimes\vabs^{\frac{1-k}{2}}\right)\times\dots\times\left(\sigma\otimes\vabs^{\frac{k-1}{2}}\right)\]
has a unique irreducible quotient denoted $\St_k(\sigma)$, which is essentially square integrable. Moreover, every irreducible, essentially square integrable representation of $\GL_n(F)$ is of the form $\St_k(\sigma)$ for a unique pair $(k,\sigma)$. Under the local Langlands correspondence $\mathscr{L}$, $\St_k(\sigma)$ corresponds to $\Sp_k\otimes\mathscr{L}(\sigma\otimes\vabs^{\frac{1-k}{2}})$ (see page 252 of \cite{ht} or Section 4.4 of \cite{rodier}), where $\Sp_k$ is as in \cite[4.1.4]{tate}. Suppose now that $n=n_1+\dots+n_r$ and that $\pi_i$ is an irreducible essentially square integrable representation of $\GL_{n_i}(F)$. Then $\pi_1\times\dots\times\pi_r$ has a distinguished constituent appearing with multiplicity one, called the Langlands subquotient, which we denote by
\[ \pi_1\boxplus\cdots\boxplus\pi_r. \]
Every irreducible representation of $\GL_n(F)$ over $\overline\Q_\ell$ is of this form for some partition of $n$, and the $\pi_i$ are well determined modulo permutation (\cite[6.1]{zelevII}). The $\pi_i$ can be ordered in such a way that the Langlands subquotient is actually a quotient of the parabolic induction.

If $\chi_1,\dots,\chi_n$ are unramified characters then 
\[ \chi_1\boxplus\cdots\boxplus\chi_n \]
is the unique unramified constituent of $\chi_1\times\dots\times\chi_n$, and every irreducible unramified representation of $\GL_n(F)$ over $\overline\Q_\ell$ is of this form. Let $\pi$ be such a representation, corresponding to a $\overline\Q_\ell$-algebra morphism $\lambda_\pi:\mathscr{H}\to\overline\Q_\ell$. For $j=1,\dots,n$, let $s_j$ denote the $j$-th elementary symmetric polynomial in $n$ variables. If we define unramified characters
\[ \chi_i:F^\times\to \overline\Q_\ell^\times \]
in such a way that $\lambda_{\pi}(T_F^{(j)})=q^{j(n-j)/2}s_j(\chi_1(\overline\omega),\dots,\chi_n(\overline\omega))$, then
\[ \pi\simeq\chi_1\boxplus\cdots\boxplus\chi_n.\]
Moreover, by the Iwasawa decomposition $\GL_n(F)=B\GL_n(\Or_F)$, we have that $\dim_{\overline\Q_\ell}\pi^{\GL_n(\Or_{F})}=1$. We denote $\mathbf{t}_\pi=(\lambda_{\pi}(T_F^{(1)}),\dots,\lambda_\pi(T_F^{(n)}))$.

\begin{lemma} Let $\pi$ be an irreducible unramified representation of $\GL_n(F)$ over $\overline\Q_\ell$. Then the characteristic polynomial of $r_\ell(\pi)(\Frob_F)$ is $P_{q,\mathbf{t}_\pi}$.
\begin{proof}
Suppose that $\pi=\chi_1\boxplus\dots\boxplus\chi_n$. Then
\[ r_\ell(\pi)=\bigoplus_{i=1}^n(\chi_i\otimes\vabs^{(1-n)/2})\circ\Art_{F}^{-1}.\]
Thus, the characteristic polynomial of $r_\ell(\pi)(\Frob_F)$ is
\[ \prod_{i=1}^n(X-\chi_i(\overline\omega)q^{(n-1)/2})=\sum_{j=0}^n(-1)^js_j(\chi_1(\overline\omega)q^{(n-1)/2},\dots,\chi_n(\overline\omega)q^{(n-1)/2})X^{n-j}=P_{q,\mathbf{t}_\pi}. \]
\end{proof}
\end{lemma}
 
Let $n=n_1+\dots +n_r$ be a partition of $n$ and let $\beta_1,\dots\beta_r$ be {\em distinct} unramified $\overline\F_\ell$-valued characters of $F^\times$. Suppose that $q\equiv 1(\Mod\ell)$. Then the representation $\beta_1[n_1]\times\dots\times\beta_r[n_r]$ is irreducible and unramified, and every irreducible unramified $\overline\F_\ell$-representation of $\GL_n(F)$ is obtained in this way. This is proved by Vigneras in \cite[VI.3]{vigbeta}. Moreover, if $\pi=\beta_1[n_1]\times\dots\times\beta_r[n_r]$, then $\pi$ is an unramified subrepresentation of the principal series $\beta_1\times\dots\times\beta_1\times\dots\times\beta_r\times\dots\times\beta_r$, where $\beta_i$ appears $n_i$ times. The Iwasawa decomposition implies that the dimension of the $\GL_n(\Or_F)$-invariants of this unramified principal series is one, and thus the same is true for $\pi$.

A character $\chi$ of $F^\times$ is called {\em tamely ramified} if it is trivial on $1+\lambda_F$, that is, if its conductor is $\leq 1$. In this case, $\chi|_{\Or_F^\times}$ has a natural extension to $U_0$, which we denote by $\chi^0$.

\begin{lemma}\label{lemainv} Let $\chi_1,\dots,\chi_n$ be $R$-valued characters of $F^\times$ such that $\chi_1,\dots,\chi_{n-1}$ are unramified and $\chi_n$ is tamely ramified. Then
\[ \dim_R\Hom_{U_0}(\chi_n^0,\chi_1\times\dots\times\chi_n)=\left\{\begin{array}{ll}n & \text{if }\chi_n\text{ is unramified}\\1 & \text{otherwise.}\end{array}\right. \]
Furthermore, if $\chi_n$ is ramified then $(\chi_1\times\dots\times\chi_n)^{U_0}=0$.
\begin{proof}
Let
\[ M(\chi_n^0)=\{f:\GL_n(\Or_{F})\to R:f(bku)=\chi(b)\chi_n^0(u)f(k)\hspace{1mm}\forall\hspace{1mm}b\in B_0,k\in\GL_n(\Or_{F}),u\in U_0\},\]
where we write $\chi$ for the character of $(F^\times)^n$ given by $\chi_1,\dots,\chi_n$.
Then, $\Hom_{U_0}(\chi_n^0,\chi_1\times\dots\times\chi_n)=(\chi_1\times\dots\times\chi_n)^{U_0=\chi_n^0}$, which by the Iwasawa decomposition is isomorphic to $M(\chi_n^0)$. By the Bruhat decomposition, 
\[ B_0\backslash\GL_n(\Or_{F})/U_0\simeq r(B_0)\backslash\GL_n(k_F)/r(U_0)\simeq W_n/W_{n-1},\]
where $W_j$ is the Weyl group of $\GL_j$ with respect to its standard maximal split torus. Here we see $W_{n-1}$ inside $W_n$ in the natural way. Let $X$ denote a set of coset representatives of $W_n/W_{n-1}$, so that
\[ \GL_n(\Or_{F})=\coprod_{w\in X}B_0wU_0.\]
Thus, if $f\in M(\chi_n^0)$, $f$ is determined by its restriction to the cosets $B_0wU_0$. We have that
\[ M(\chi_n^0)\simeq\prod_{w\in X}M_w,\]
where $M_w$ is the space of functions on $B_0wU_0$ satisfying the transformation rule of $M(\chi_n^0)$. It is clear that $\dim_RM_w\leq 1$ for every $w$. Moreover, if $\chi_n$ is unramified, then $M_w$ is non-zero, a non-zero function being given by $f(w)=1$. Thus, in this case, $\dim_RM(\chi_n^0)=n$. 

In the ramified case, let $a=\diag(a_1,\dots,a_n)\in B_0$, with $a_i\in\Or_F^\times$ and $a_n$ such that $\chi_n(a_n)\neq 1$. Then 
\[ \chi_n(a_n)f(w)=f(aw)=f(wa^w)=\chi_n^0(a^w)f(w)=f(w) \]
unless $w\in W_{n-1}$. Thus, only the identity coset survives, and $\dim_RM(\chi_n^0)=1$.

For the last assertion, let $f\in(\chi_1\times\dots\times\chi_n)$ be $U_0$-invariant. To see that it is zero, it is enough to see that $f(w)=0$ for every $w\in X$. Choosing $a\in\GL_n(\Or_F)$ to be a scalar matrix corresponding to an element $a\in\Or_F^\times$ for which $\chi_n(a)\neq 1$, we see that $a$ is in $B_0$ (and hence in $U_0$), thus $f(aw)=\chi_n(a)f(w)=f(wa)=f(w)$, so $f(w)=0$ for any $w\in X$.
\end{proof}
\end{lemma}

Let $P_M$ denote the parabolic subgroup of $\GL_n(F)$ containing $B$ corresponding to the partition $n=(n-1)+1$, and let $U_M$ denote its unipotent radical. Take the Levi decomposition $P_M=MU_M$, where $M\simeq\GL_{n-1}(F)\times\GL_1(F)$. Consider the opposite parabolic subgroup $\overline{P_M}$ with Levi decomposition $\overline{P_M}=M\overline{U_M}$. Let
\[ U_{0,M}=U_0\cap M\simeq\GL_{n-1}(\Or_F)\times\GL_1(\Or_F).\]
Let $\chi_n$ be a tamely ramified character of $F^\times$, and let $\chi_n^0$ be its extension to $U_0$. Let
\[ \mathscr{H}_M(\chi_n)=\End_{M}(\ind_{U_{0,M}}^M\chi_n),\]
where $\ind$ denotes compact induction and $\chi_n$ is viewed as a character of $U_{0,M}$ via projection to the last element of the diagonal. Thus, $\mathscr{H}_M(\chi_n)$ can be identified with the $R$-vector space of compactly supported functions $f:M\to R$ such that $f(kmk')=\chi_n(k)f(m)\chi_n(k')$ for $m\in M$ and $k,k'\in U_{0,M}$. Similarly, let
\[ \mathscr{H}_0(\chi_n)=\End_{\GL_n(F)}(\ind_{U_0}^{\GL_n(F)}\chi_n^0).\]
This is identified with the $R$-vector space of compactly supported functions $f:\GL_n(F)\to R$ such that $f(kgk')=\chi_n^0(k)f(g)\chi_n^0(k')$ for every $g\in\GL_n(F)$, $k,k'\in U_0$. There is a natural injective homomorphism of $R$-modules
\[ \mathscr{T}:\mathscr{H}_M(\chi_n)\to\mathscr{H}_0(\chi_n),\]
which can be described as follows (see \cite[II.3]{vigbeta}). Let $m\in M$. Then $\mathscr{T}(1_{U_{0,M}mU_{0,M}})=1_{U_0mU_0}$, where $1_{U_{0,M}mU_{0,M}}$ is the function supported in $U_{0,M}mU_{0,M}$ whose value at $umu'$ is $\chi_n(u)\chi_n(u')$, and similarly for $1_{U_0mU_0}$. Define
\[ U_0^+=U_0\cap U_M \]
and
\[ U_0^-=U_0\cap\overline{U_M}.\]
Then $U_0=U_0^-U_{0,M}U_0^+=U_0^+U_{0,M}U_0^-$, and $\chi_n^0$ is trivial on $U_0^-$ and $U_0^+$.
Let
\[ M^-=\{m\in M/\hspace{1mm}m^{-1}U_0^+m\subset U_0^+\hspace{1mm}\text{ and }\hspace{1mm}mU_0^-m^{-1}\subset U_0^-\}.\]
We denote by $\mathscr{H}_M^-(\chi_n)$ the subspace of $\mathscr{H}_M(\chi_n)$ consisting of functions supported on the union of cosets of the form $U_{0,M}mU_{0,M}$ with $m\in M^-$. 

\begin{prop}\label{subalgplus}
The subspace $\mathscr{H}_M^-(\chi_n)\subset\mathscr{H}_M(\chi_n)$ is a subalgebra, and the restriction $\mathscr{T}^-:\mathscr{H}_M^-(\chi_n)\to\mathscr{H}_0(\chi_n)$ is an algebra homomorphism.
\begin{proof} This is proved in \cite[II.5]{vigbeta}.
\end{proof}
\end{prop}

Let $\pi$ be a representation of $\GL_n(F)$ over $R$. Then $\Hom_{\GL_n(F)}(\ind_{U_0}^{\GL_n(F)}\chi_n^0,\pi)$ is naturally a right module over $\mathscr{H}_0(\chi_n)$. By the adjointness between compact induction and restriction,
\[ \Hom_{\GL_n(F)}(\ind_{U_0}^{\GL_n(F)}\chi_n^0,\pi)=\Hom_{U_0}(\chi_n^0,\pi),\]
and therefore the right hand side is also a right $\mathscr{H}_0(\chi_n)$-module. There is an $R$-algebra isomorphism $\mathscr{H}_0(\chi_n)\simeq\mathscr{H}_0(\chi_n^{-1})^{\opp}$ given by $f\mapsto f^*$, where $f^*(g)=f(g^{-1})$. We then see $\Hom_{U_0}(\chi_n^0,\pi)$ as a left $\mathscr{H}_0(\chi_n^{-1})$-module in this way. Similarly, $\Hom_{U_{0,M}}(\chi_n,\pi)$ is a left $\mathscr{H}_M(\chi_n^{-1})$-module when $\pi$ is a representation of $M$ over $R$. For a representation $\pi$ of $\GL_n(F)$, let $\pi_{\overline{U_M}}$ be the representation of $M$ obtained by (non-normalized) parabolic restriction. Then the natural projection $\pi\to\pi_{\overline{U_M}}$ is $M$-linear.

\begin{rem} Let $\overline{B_{n-1}}$ denote the subgroup of lower triangular matrices of $\GL_{n-1}(F)$, so that $\overline{B_{n-1}}\times\GL_1(F)$ is a parabolic subgroup of $M$, with the standard maximal torus $T\subset M$ of $\GL_n(F)$ as a Levi factor. Let $\chi_1,\dots,\chi_n$ be characters of $F^\times$. Then
\begin{equation}
\label{eqss}
\left((\chi_1\times\dots\times\chi_n)_{\overline{U_M}}\right)^{\ssp}\simeq\bigoplus_{i=1}^n\left(i_{\overline{B_{n-1}}\times\GL_1(F)}^M(\chi^{w_i})\right)^{\ssp}\otimes\delta_{\overline{P_M}}^{1/2},
\end{equation}
where $\ssp$ denotes semisimplification and $i_{\overline{B_{n-1}}\times\GL_1(F)}^M$ is the normalized parabolic induction. Here, $w_i$ is the permutation of $n$ letters such that $w_i(n)=n+1-i$ and $w_i(1)>w_i(2)>\dots>w_i(n-1)$. This follows from Theorem 6.3.5 of \cite{casselman} when $R=\overline\Q_\ell$. As Vign\'eras points out in \cite[II.2.18]{vigbeta}, the same proof is valid for the $R=\overline\F_\ell$ case.\end{rem}

\begin{prop}\label{proj} Let $\chi_1,\dots,\chi_n$ be $R$-valued characters of $F^\times$, such that $\chi_1,\dots,\chi_{n-1}$ are unramified and $\chi_n$ is tamely ramified.
\begin{enumerate}[(i)]
\item The natural projection $\chi_1\times\dots\times\chi_n\to(\chi_1\times\dots\times\chi_n)_{\overline{U_M}}$ induces an isomorphism of $R$-modules
 \begin{equation}
 \label{eqp}
  p:\Hom_{U_0}(\chi_n^0,(\chi_1\times\dots\times\chi_n))\to\Hom_{U_{0,M}}(\chi_n,(\chi_1\times\dots\times\chi_n)_{\overline{U_M}}).
 \end{equation}
\item For every $\phi\in\Hom_{U_0}(\chi_n^0,(\chi_1\dots\times\dots\chi_n))$ and every $m\in M^-$,
\[ p(1_{U_0mU_0}.\phi)=\delta_{P_M}(m)1_{U_{0,M}mU_{0,M}}.p(\phi).\]
\end{enumerate}
\begin{proof} The last assertion is proved in \cite[II.9]{vigbeta}. The fact that $p$ is surjective follows by \cite[II.3.5]{vigirr}. We prove injectivity now. By Lemma \ref{lemainv}, the dimension of the left hand side is $n$ if $\chi_n$ is unramified and $1$ otherwise. Suppose first that $R=\overline\Q_\ell$. If $\chi_n$ is unramified, each summand of the right hand side of (\ref{eqss}) has a one dimensional $U_{0,M}$-fixed subspace, while if $\chi_n$ is ramified, only the summand corresponding to the identity permutation has a one dimensional $U_{0,M}$-fixed subspace, all the rest being zero. This implies that
\[ \dim_{\overline\Q_\ell}\left((\chi_1\times\dots\times\chi_n)_{\overline{U_M}}\right)^{U_{0,M}}=\left\{\begin{array}{ll}n & \text{if }\chi_n\text{ is unramified}\\1 & \text{otherwise,}\end{array}\right.\]
Therefore $p$ is an isomorphism for reasons of dimension. This completes the proof of the injectivity of $p$ over $\overline\Q_\ell$.

We give the proof over $\overline\F_\ell$ only in the unramified case, the ramified case being similar. First of all, note that the result for $\overline\Q_\ell$ implies the corresponding result over $\overline\Z_\ell$, the ring of integers of $\overline\Q_\ell$. Indeed, suppose each $\chi_i$ takes values in $\overline\Z_\ell^\times$, and let $(\chi_1\times\dots\times\chi_n)_{\overline\Z_\ell}$ (respectively, $(\chi_1\times\dots\times\chi_n)_{\overline\Q_\ell}$) denote the parabolic induction over $\overline\Z_\ell$ (respectively, $\overline\Q_\ell$). Then $(\chi_1\times\dots\times\chi_n)_{\overline\Z_\ell}$ is a lattice in $(\chi_1\times\dots\times\chi_n)_{\overline\Q_\ell}$, that is, a free $\overline\Z_\ell$-submodule which generates $(\chi_1\times\dots\times\chi_n)_{\overline\Q_\ell}$ and is $\GL_n(F)$-stable (\cite[II.4.14(c)]{vigirr}). It then follows that $((\chi_1\times\dots\times\chi_n)_{\overline\Z_\ell})^{U_0}$ is a lattice in $(\chi_1\times\dots\times\chi_n)_{\overline\Q_\ell})^{U_0}$ (\cite[I.9.1]{vigirr}), and so is free of rank $n$ over $\overline\Z_\ell$. Simiarly, $((\chi_1\times\dots\times\chi_n)_{\overline{U_M},\overline\Z_\ell})^{U_{0,M}}$ is a lattice in $((\chi_1\times\dots\times\chi_n)_{\overline{U_M},\overline\Q_\ell})^{U_{0,M}}$ (\cite[II.4.14(d)]{vigirr}), and thus it is free of rank $n$ over $\overline\Z_\ell$. Moreover, the map $p$ with coefficients in $\overline\Z_\ell$ is still surjective (\cite[II 3.3]{vigirr}), hence it is an isomorphism by reasons of rank.

Finally, consider the $\overline\F_\ell$ case. Choose liftings $\widetilde\chi_i$ of $\chi_i$ to $\overline\Z_\ell$-valued characters. Then there is a natural injection
\[
(\widetilde\chi_1\times\dots\times\widetilde\chi_n)_{\overline{U_M}}\otimes_{\overline\Z_\ell}\overline\F_\ell\hookrightarrow (\chi_1\times\dots\times\chi_n)_{\overline{U_M}}
\]
inducing an injection
\begin{equation}\label{eqpmodl}
((\widetilde\chi_1\times\dots\times\widetilde\chi_n)_{\overline{U_M}})^{U_{0,M}}\otimes_{\overline\Z_\ell}\overline\F_\ell\hookrightarrow ((\chi_1\times\dots\times\chi_n)_{\overline{U_M}})^{U_{0,M}}.
\end{equation}
Now, we have seen that the left hand side of (\ref{eqpmodl}) has dimension $n$ over $\overline\F_\ell$. We claim that the right hand side of (\ref{eqpmodl}) has dimension $\leq n$. Indeed, by looking at the right hand side of (\ref{eqss}), this follows from the fact that the $U_{0,M}$-invariants of the semisimplification have dimension $n$. Thus, (\ref{eqpmodl}) is an isomorphism and $\dim_{\overline\F_\ell}(\chi_1\times\dots\times\chi_n)_{\overline{U_M}})^{U_{0,M}}=n$. Since the left hand side of (\ref{eqp}) has dimension $n$ and $p$ is surjective, it must be an isomorphism.
\end{proof}
\end{prop}

Let $\mathscr{H}_0$ (respectively, $\mathscr{H}_1$) be the $R$-valued Hecke algebra of $\GL_n(F)$ with respect to $U_0$ (respectively, $U_1$). Thus, $\mathscr{H}_0=\mathscr{H}_0(1)$. If $\pi$ is a representation of $\GL_n(F)$ over $R$, then $\pi^{U_0}$ is naturally a left $\mathscr{H}_0$-module. For any $\alpha\in F^\times$ with $|\alpha|\leq 1$, let $m_{\alpha}\in M$ be the element
\[ m_{\alpha}=\left(\begin{array}{cc}1_{n-1}&0\\0&\alpha\end{array}\right).\]
For $i=0$ or $1$, let $V_{\alpha,i}\in\mathscr{H}_i$ be the Hecke operators $\left[U_im_{\alpha}U_i\right]$. If $\pi$ is a representation of $\GL_n(F)$, then $\pi^{U_0}\subset\pi^{U_1}$ and the action of the operators defined above are compatible with this inclusion.

Let $\mathscr{H}_M=\mathscr{H}_M(1)$, and let $V_{\overline\omega,M}=\left[U_{0,M}m_{\overline\omega}U_{0,M}\right]\in\mathscr{H}_M$. Since $m_{\overline\omega}\in M^-$, $V_{\overline\omega,M}\in\mathscr{H}_M^-$, and $\mathscr{T}^-(V_{\overline\omega,M})=V_{\overline\omega,0}\in\mathscr{H}_0$. As above, if $\pi$ is a representation of $M$ over $R$, we consider the natural left action $\mathscr{H}_M$ on $\pi^{U_{0,M}}$.

\begin{coro} Let $\chi_1,\dots,\chi_n$ be $\overline\Q_\ell$-valued unramified characters of $F^\times$. Then the set of eigenvalues of $V_{\overline\omega,0}$ acting on the $n$-dimensional space $(\chi_1\times\dots\times\chi_n)^{U_0}$ is equal (counting multiplicities) to $\{q^{(n-1)/2}\chi_i(\overline\omega)\}_{i=1}^n$.
\begin{proof} Note that $V_{\overline\omega,M}$ acts on the $U_{0,M}$-invariants of each summand of the right hand side of (\ref{eqss}) by the scalar $\chi_i(\overline\omega)q^{(1-n)/2}$. Thus, the eigenvalues of $V_{\overline\omega,M}$ in $(\chi_1\times\dots\times\chi_n)_{\overline{U_M}}^{U_{0,M}}$ are the $q^{(1-n)/2}\chi_i(\overline\omega)$. The corollary follows then by Proposition \ref{proj}.
\end{proof}
\end{coro}

\begin{prop}
\label{eigenV} Let $\pi$ be an irreducible unramified representation of $\GL_n(F)$ over $R$. Then $\pi^{U_0}=\pi^{U_1}$ and the following properties hold.
\begin{enumerate}[(i)]
\item If $R=\overline\Q_\ell$ and $\pi=\chi_1\boxplus\dots\boxplus\chi_n$, with $\chi_i$ unramified characters of $F^\times$, then $\dim_R\pi^{U_0}\leq n$ and the eigenvalues of $V_{\overline\omega,0}$ acting on $\pi^{U_0}$ are contained in $\{q^{(n-1)/2}\chi_i(\overline\omega)\}_{i=1}^n$ (counting multiplicities).
\item If $R=\overline\F_\ell$, $q\equiv1(\Mod\ell)$ and $\pi=\beta_1[n_1]\times\dots\times\beta_r[n_r]$ with $\beta_i$ distinct unramified characters of $F^\times$, then $\dim_R\pi^{U_0}=r$ and $V_{\overline\omega,0}$ acting on $\pi^{U_0}$ has the $r$ distinct eigenvalues $\{\beta_j(\overline\omega)\}_{j=1}^r$.
\end{enumerate}
\begin{proof}
The fact that $\pi^{U_1}=\pi^{U_0}$ follows immediately because the central character of $\pi$ is unramified. Since taking $U_0$-invariants is exact in characteristic zero, part (i) is clear from the last corollary. Let us prove (ii). Let $P$ be the parabolic subgroup of $\GL_n(F)$ containing $B$ corresponding to the partition $n=n_1+\cdots+n_r$. As usual, since $\GL_n(F)=P\GL_n(\Or_F)$, the $\overline\F_\ell$-dimension of $\pi^{U_0}$ is equal to the cardinality of $\left(\GL_n(\Or_F)\cap P\right)\backslash\GL_n(\Or_F)/U_0$. By the Bruhat decomposition, this equals the cardinality of
\[ \mathfrak{S}_{n_1}\times\cdots\times\mathfrak{S}_{n_r}\backslash\mathfrak{S}_n/\mathfrak{S}_{n-1}\times\mathfrak{S}_1,\]
where $\mathfrak{S}_i$ is the symmetric group on $i$ letters. This cardinality is easily seen to be $r$.

It remains to prove the assertion about the eigenvalues of $V_{\overline\omega,0}$ on $\pi^{U_0}$. Let us first replace $U_0$ by $\Iw$ (this was first suggested by Vign\'eras). By the Iwasawa decomposition and the Bruhat decomposition,
\[ \GL_n(F)=\coprod_{s\in S}Ps\Iw,\]
where $S\subset\GL_n(F)$ is a set of representatives for $(\mathfrak{S}_{n_1}\times\dots\times\mathfrak{S}_{n_r})\backslash\mathfrak{S}_n$. Then $\pi^{\Iw}$ has as a basis the set $\{\varphi_s\}_{s\in S}$, where $\varphi_s$ is supported on $Ps\Iw$ and $\varphi_s(s)=1$. 

Let $\mathscr{H}_{\overline\F_\ell}(n,1)$ denote the Iwahori-Hecke algebra for $\GL_n(F)$ over $\overline\F_\ell$, that is, the Hecke algebra for $\GL_n(F)$ with respect to the compact open subgroup $\Iw$. Thus, $\pi^{\Iw}$ is naturally a left module over $\mathscr{H}_{\overline\F_\ell}(n,1)$. For $i=1,\dots,n-1$, let $s_i$ denote the $n$ by $n$ permutation matrix corresponding to the transposition $(i\hspace{2mm}i+1)$, and let $S_i=[\Iw s_i\Iw]\in\mathscr{H}_{\overline\F_\ell}(n,1)$. For $j=0,\dots,n$, let $t_j$ denote the diagonal matrix whose first $j$ coordinates are equal to $\overline\omega$, and whose last $n-j$ coordinates are equal to $1$. Let $T_j=[\Iw t_j\Iw]\in\mathscr{H}_{\overline\F_\ell}(n,1)$, and for $j=1,\dots,n$, let $X_j=T_j(T_{j-1}^{-1})$. Then $\mathscr{H}_{\overline\F_\ell}(n,1)$ is generated as an $\overline\F_\ell$-algebra by $\{S_i\}_{i=1}^{n-1}\cup\{X_1,X_1^{-1}\}$ (\cite[I.3.14]{vigirr}). We denote by $\mathscr{H}_{\overline\F_\ell}^0(n,1)$ the subalgebra generated by $\{S_i\}_{i=1}^{n-1}$, which is canonically isomorphic to the group algebra $\overline\F_\ell[\mathfrak{S}_n]$ of the symmetric group (\cite[I.3.12]{vigirr}). It can also be identified with the Hecke algebra of $\GL_n(\Or_F)$ with respect to $\Iw$ (\cite[I.3.14]{vigirr}). The subalgebra $A=\overline\F_\ell[\{X_i^\pm\}_{i=1}^n]$ is commutative, and characters of $T$ can be seen as characters on $A$. Let $\chi_1,\dots,\chi_n:F^\times\to\overline\F_\ell^\times$ be the characters defined by
\[
\chi_1=\cdots=\chi_{n_1}=\beta_1;\]
\[\cdots;\]
\[\chi_{n_1+\cdots+n_{j-1}+1}=\cdots=\chi_{n_1+\cdots+n_j}=\beta_j;\]
\[\cdots.\]
Then the action of $A$ on $\varphi_s$ is given by the character $s(\chi)$. Note that the set $\{s(\chi)\}_{s\in S}$ is just the set of $n$-tuples of characters in which $\beta_i$ occurs $n_i$ times, with arbitrary order. It is clear that for each $j=1,\dots,r$, there is at least one $s\in S$ for which $s(n)\in\{n_1+\dots+n_{j-1}+1,\dots,n_1+\dots+n_j\}$, so that $X_n\varphi_s=\beta_j(\overline\omega)\varphi_s$. Let
\[ \varphi=\sum_{s\in S}\varphi_s.\]
Then $\varphi$ generates $\pi^{\GL_n(\Or_F)}$. For $j=1,\dots,r$, let
\[ \psi_j=\sum_{s\in S,\chi_{s(n)}=\beta_j}\varphi_s.\]
We have seen above that $\psi_j\neq 0$. Moreover, $X_n\psi_j=\beta_j(\overline\omega)\psi_j$. Let $P_j\in\overline\F_\ell[X]$ be a polynomial such that $P_j(\beta_j(\overline\omega))=1$ and $P_j(\beta_i(\overline\omega))=0$ for every $i\neq j$. Then $\psi_j=P_j(X_n)\varphi$, and it follows that the $r$ distinct eigenvalues $\{\beta_j(\overline\omega)\}_{j=1}^r$ of $X_n$ on $\pi^{\Iw}$ already occur on the subspace $\overline\F_\ell[X_n]\varphi$.

Consider now the map $p_T:\pi^{\Iw}\to(\pi_{\overline N})^{T_0}$, where $\overline N$ is the unipotent radical of the parabolic subgroup of $\GL_n(F)$ containing $T$, opposite to $B$, and $T_0=T\cap\GL_n(\Or_F)$. By \cite[II.3.5]{vigirr}, $p_T$ is an isomorphism. On the other hand, there is a commutative diagram
\[\xymatrix{
 \pi^{U_0}\ar[r]^i\ar[d]_{p_M} & \pi^{\Iw}\ar[d]_{p_T} \\ (\pi_{\overline{U_M}})^{U_{0,M}} \ar[r]^{p_{M,T}} & (\pi_{\overline N})^{T_0},}
\]
where $i$ is the inclusion and $p_M$ and $p_{M,T}$ are the natural projection to the coinvariants. The analogues of part (ii) of Proposition \ref{proj} for $p_M$, $p_T$ and $p_{M,T}$ are still valid (\cite[II.9]{vigbeta}). Thus, for $f\in\pi^{U_0}$,
\[ p_{T}(i(V_{\overline\omega,0}f))=p_{M,T}(p_M(V_{\overline\omega,0}f))=p_{M,T}([U_{0,M}m_{\overline\omega}U_{0,M}]p_M(f))=\]
\[=[T_0m_{\overline\omega}T_0]p_{M,T}(p_M(f))=[T_0m_{\overline\omega}T_0]p_T(i(f))=p_T(X_ni(f)).\]
It follows that $V_{\overline\omega,0}=X_n$ on $\pi^{U_0}$. In particular, $\overline\F_\ell[X_n]\varphi=\overline\F_\ell[V_{\overline\omega,0}]\varphi\subset\pi^{U_0}$. By what we have seen above, we conclude that the eigenvalues of $V_{\overline\omega,0}$ on the $r$ dimensional space $\pi^{U_0}$ are $\{\beta_j(\overline\omega)\}_{j=1}^r$, as claimed.
\end{proof}
\end{prop}

\begin{coro}\label{polV} Suppose that $q\equiv1(\Mod\ell)$ and let $\pi$ be an irreducible unramified representation of $\GL_n(F)$ over $\overline\F_\ell$. Let $\varphi\in\pi^{\GL_n(\Or_F)}$ be a non-zero spherical vector. Then $\varphi$ generates $\pi^{U_0}$ as a module over the algebra $\overline\F_\ell[V_{\overline\omega,0}]$.
\begin{proof} This is actually a corollary of the proof of the above proposition. Indeed, $V_{\overline\omega,0}$ has $r$ distinct eigenvalues on $\overline\F_\ell[V_{\overline\omega,0}]\varphi\subset\pi^{U_0}$, and $\dim_{\overline\F_\ell}\pi^{U_0}=r$.
\end{proof}
\end{coro}

\begin{lemma}\label{lemma315} Let $\pi$ be an irreducible representation of $\GL_n(F)$ over $\overline\Q_\ell$ with a non-zero $U_1$-fixed vector but no non-zero $\GL_n(\Or_F)$-fixed vectors. Then $\dim_{\overline\Q_\ell}\pi^{U_1}=1$ and there is a character
 \[ V_\pi:F^\times\to\overline\Q_\ell^\times \]
with open kernel such that for every $\alpha\in F^\times$ with non-negative valuation, $V_\pi(\alpha)$ is the eigenvalue of $V_{\alpha,1}$ on $\pi^{U_1}$.  Moreover, there is an exact sequence
\[ 0\longrightarrow s\longrightarrow r_\ell(\pi)\longrightarrow V_\pi\circ\Art_{F}^{-1}\longrightarrow 0,\]
where $s$ is unramified. If $\pi^{U_0}\neq 0$ then $q^{-1}V_\pi(\overline\omega)$ is a root of the characteristic polynomial of $s(\Frob_F)$. If, on the other hand, if $\pi^{U_0}=0$, then $r_\ell(\pi)(\Gal(\overline F/F))$ is abelian.

\begin{proof} This is Lemma 3.1.5 of \cite{cht}. The proof basically consists in noting that if $\pi^{U_1}\neq 0$, then either $\pi\simeq\chi_1\boxplus\dots\boxplus\chi_n$ with $\chi_1,\dots,\chi_{n-1}$ unramified and $\chi_n$ tamely ramified, or $\pi\simeq\chi_1\boxplus\dots\boxplus\chi_{n-2}\boxplus\St_2(\chi_{n-1})$ with $\chi_1,\dots,\chi_{n-1}$ unramified. Then one just analyzes the cases separately, and calculates explicitly the action of the operators $U_{F,1}^{(j)}$ (see \cite{cht} for their definition) and $V_{\alpha,1}$.
\end{proof}

\end{lemma}

\begin{lemma}\label{321} Suppose that $q\equiv1(\Mod\ell)$, and let $\pi$ be an irreducible unramified representation of $\GL_n(F)$ over $\overline\F_\ell$. Let $\lambda_\pi(T_F^{(j)})$ be the eigenvalue of $T_F^{(j)}$ on $\pi^{\GL_n(\Or_F)}$, and $\mathbf{t}_\pi=(\lambda_\pi(T_F^{(1)}),\dots,\lambda_\pi(T_F^{(n)}))$. Suppose that $P_{q,\mathbf{t}_\pi}=(X-a)^mF(X)$ in $\overline\F_\ell[X]$, with $m>0$ and $F(a)\neq 0$. Then $F(V_{\overline\omega,0})$, as an operator acting on $\pi^{U_0}$, is non-zero on the subspace $\pi^{\GL_n(\Or_F)}$.
\begin{proof} Suppose on the contrary that $F(V_{\overline\omega,0})(\pi^{\GL_n(\Or_F)})=0$. Let $\varphi\in\pi^{\GL_n(\Or_F)}$ be a non-zero element. Suppose that $\pi=\beta_1[n_1]\times\dots\times\beta_r[n_r]$, with $\beta_i$ distinct unramified $\overline\F_\ell^\times$-valued characters of $F^\times$. Then, since $q=1$ in $\overline\F_\ell$,
\[ P_{q,\mathbf{t}_\pi}=\prod_{i=1}^r(X-\beta_i(\overline\omega))^{n_i}.\]
Suppose that $a=\beta_j(\overline\omega_1)$, so that $F(X)=\prod_{i\neq j}(X-\beta_i(\overline\omega))^{n_i}$. By Proposition \ref{eigenV} (ii), $\pi^{U_0}$ has dimension $r$ and $V_{\overline\omega,0}$ is diagonalizable on this space, with distinct eigenvalues $\beta_i(\overline\omega)$. Let $\varphi_j\in\pi^{U_0}$ denote an eigenfunction of $V_{\overline\omega,0}$ of eigenvalue $\beta_j(\overline\omega)$. By Corollary \ref{polV}, there exists a polynomial $P_j\in\overline\F_\ell[X]$ such that $\varphi_j=P_j(V_{\overline\omega,0})(\varphi)$. Since polynomials in $V_{\overline\omega,0}$ commute with each other, we must have $F(V_{\overline\omega,0})(\varphi_j)=0$, but this also equals $F(\beta_j(\overline\omega))\varphi_j\neq 0$, which is a contradiction.
\end{proof}
\end{lemma}

\section{Automorphic forms on unitary groups}

\subsection{Totally definite groups} Let $F^+$ be a totally real field and $F$ a totally imaginary quadratic extension of $F^+$. Denote by $c\in\Gal(F/F^+)$
the non-trivial Galois automorphism. Let $n\geq 1$ be an integer and $V$ an $n$-dimensional vector space over $F$, equipped
with a non-degenerate $c$-hermitian form $h:V\times V\to F$. To the pair $(V,h)$ there is attached a reductive algebraic group $U(V,h)$
over $F^+$, whose points in an $F^+$-algebra $R$ are
\[ U(V,h)(R)=\{g\in\Aut_{(F\otimes_{F^+}R)-\lin}(V\otimes_{F+}R):h(gx,gy)=h(x,y)\hspace{1mm}\forall x,y\in V\otimes_{F^+}R\}.\]
By an {\em unitary group} attached to $F/F^+$ in $n$ variables,
we shall mean an algebraic group of the form $U(V,h)$ for some pair $(V,h)$ as above. Let $G$ be such a group. Then $G_F=G\otimes_{F^+}F$ is isomorphic to $\GL_V$, and in fact it is an outer form of $\GL_V$. Let $G(F^+_\infty)=\prod_{v|\infty}G(F^+_v)$, and if $v$ is any place of $F^+$, let $G_v=G\otimes_{F^+}F^+_v$. We say that $G$ is {\em totally definite} if $G(F^+_\infty)$ is compact (and thus isomorphic to a product of copies of the compact unitary group $U(n)$).

Suppose that $v$ is a place of $F^+$ which splits in $F$, and let $w$ be a place of $F$ above $v$, corresponding to an $F^+$-embedding $\sigma_w:F\hookrightarrow\overline{F^+_v}$. Then $F^+_v=\sigma_w(F)F^+_v$ is an $F$-algebra by means of $\sigma_w$, and thus $G_v$ is isomorphic to $\GL_{V\otimes F^+_v}$, the tensor product being over $F$. Note that if we choose another place $w^c$ of $F$ above $v$, then $\sigma_w$ and $\sigma_{w^c}$ give $F^+_v$ two different $F$-algebra structures. If we choose a basis of $V$, we obtain two isomorphisms $i_w,i_{w^c}:G_v\to\GL_{n/F^+_v}$. If $X\in\GL_n(F)$ is the matrix of $h$ in the chosen basis, then for any $F^+_v$-algebra $R$ and any $g\in G_v(R)$, $i_{w^c}(g)=X^{-1}(^ti_w(g)^{-1})X$, where we see $X\in\GL_n(R)$ via $\sigma_w:F\to F^+_v\to R$. 

The choice of a lattice $L$ in $V$ such that $h(L\times L)\subset\Or_F$ gives an affine group scheme over $\Or_{F^+}$, still denoted by $G$, which is isomorphic to $G$ after extending scalars to $F^+$. We will fix from now on a basis for $L$ over $\Or_F$, giving also
an $F$-basis for $V$; with respect to these, for each split place $v$ of $F^+$ and each place $w$ of $F$ above $v$, $i_w$ gives an isomorphism between $G(F^+_v)$ and
$\GL_n(F_w)$ taking $G(\Or_{F^+_v})$ to $\GL_n(\Or_{F_w})$.

\subsection{Automorphic forms} Let $G$ be a totally definite unitary group in $n$ variables attached
to $F/F^+$. We let $\mathscr{A}$ denote the space of automorphic forms on $G(\A_{F^+})$. Since the group is totally definite, $\mathscr{A}$
decomposes, as a representation of $G(\A_{F^+})$, as
\[ \mathscr{A}\cong\bigoplus_{\pi}m(\pi)\pi,\]
where $\pi$ runs through the isomorphism classes of irreducible admissible representations of $G(\A_{F^+})$, and $m(\pi)$ is the multiplicity
of $\pi$ in $\mathscr{A}$, which is always finite. This is a well known fact for any reductive group compact at infinity, but we recall the proof as a warm up for the following sections and to set some notation. The isomorphism classes of continuous, complex, irreducible (and hence finite dimensional) representations of $G(F^+_\infty)$ are parametrized by elements $\bbc=(b_\tau)\in(\Z^{n,+})^{\Hom(F^+,\R)}$. We denote them by $W_\bbc$. Since $G(F^+_\infty)$ is compact and every element of $\mathscr{A}$ is $G(F^+_\infty)$-finite, $\mathscr{A}$ decomposes as a direct sum of irreducible $G(\A_{F^{+}})$-submodules. Moreover, we can write
\[ \mathscr{A}\cong\bigoplus_{\bbc}W_\bbc\otimes\Hom_{G(F^+_\infty)}(W_\bbc,\mathscr{A})\]
as $G(\A_{F^+})$-modules. Denote by $\A_{F^+}^\infty$ the ring of finite ad{\`e}les. For any $\bbc$, let $S_{\bbc}$ be the space of smooth (that is, locally constant) functions $f:G(\A_{F^+}^\infty)\to W_\bbc^\vee$ such that
$f(\gamma g)=\gamma_\infty f(g)$ for all $g\in G(\A_{F^+}^\infty)$ and $\gamma\in G(F^+)$. Then the map
\[f\mapsto\left(w\mapsto\left(g\mapsto(g_\infty^{-1}f(g^\infty))(w)\right)\right) \]
induces a $G(\A_{F^+}^\infty)$-isomorphism between $\Hom_{G(F^+_\infty)}(W_\bbc,\mathscr{A})$ and $S_{\bbc}$, where the action on this last space is
by right translation. For every compact open subgroup $U\subset G(\A_{F^+}^\infty)$, the space $G(F)\backslash G(\A_{F^+}^\infty)/U$ is finite,
and hence the space of $U$-invariants of $S_{\bbc}$ is finite-dimensional. In particular, every irreducible summand of $W_\bbc\otimes\Hom_{G(F^+_\infty)}(W_\bbc,\mathscr{A})$ is admissible and appears with finite multiplicity. Thus, every irreducible summand of $\mathscr{A}$ is admissible, and appears with finite
multiplicity because its isotypic component is contained in $W_\bbc\otimes\Hom_{G(F^+_\infty)}(W_\bbc,\mathscr{A})$ for some $\bbc$.

\subsection{$\ell$-adic models of automorphic forms}\label{sell} Let $\ell$ be an odd prime number. We will assume, from now on to the end of
this section, that every place of $F^+$ above $\ell$ splits in $F$. Let $K$ be a finite extension of $\Q_\ell$. Fix an algebraic closure $\overline K$ of $K$, and suppose that $K$ is big enough to contain all embeddings of $F$ into $\overline K$. Let $\Or$ be its ring
 of integers and $\lambda$ its maximal ideal. Let $S_\ell$ denote the set of places of $F^+$ above $\ell$, and $I_\ell$ the set of
 embeddings $F^+\hookrightarrow K$. Thus, there is a natural surjection $h:I_\ell\twoheadrightarrow S_\ell$. Let $\widetilde S_\ell$
 denote a set of places of $F$ such that $\widetilde S_\ell\coprod \widetilde S_\ell^c$ consists of all the places above $S_\ell$;
 thus, there is a bijection $S_\ell\simeq\widetilde S_\ell$. For $v\in S_\ell$, we denote by $\widetilde v$ the corresponding place in $\widetilde S_\ell$. Also, let $\widetilde I_\ell$ denote the set of embeddings $F\hookrightarrow
 K$ giving rise to a place in $\widetilde S_\ell$. Thus, there is a bijection between $I_\ell$ and $\widetilde{I_\ell}$, which we denote by $\tau\mapsto\widetilde\tau$. Also, denote by $\tau\mapsto w_\tau$ the natural surjection $\widetilde I_\ell\to\widetilde S_\ell$. Finally, Let $F^+_\ell=\prod_{v\mid\ell}F^+_v$. 

Let $\ba\in(\Z^{n,+})^{\Hom(F,K)}$. Consider the following representation of $G(F^+_\ell)\simeq\prod_{\widetilde v\in
\widetilde S_\ell}\GL_n(F_{\widetilde v})$. For each $\widetilde\tau\in\widetilde I_\ell$, we have an embedding $\GL_n(F_{w_{\widetilde\tau}})\hookrightarrow\GL_n(K)$. Taking the product over $\widetilde\tau$ and composing with the projection on the $w_{\widetilde\tau}$-coordinates, we have an irreducible representation \[ \xi_\ba:G(F^+_\ell)\longrightarrow\GL(W_{\ba}),\] where
$W_{\ba}=\otimes_{\widetilde\tau\in\widetilde I_\ell}W_{a_{\widetilde\tau},K}$. This representation has an integral model
$\xi_\ba:G(\Or_{F^+_\ell})\to\GL(M_{\ba})$. In order to base change to automorphic representations of $\GL_n$, we need to impose the additional assumption that
\[ a_{\tau c,i}=-a_{\tau,n+1-i} \]
for every $\tau\in\Hom(F,K)$ and every $i=1,\dots,n$.

Besides the weight, we will have to introduce another collection of data, away from $\ell$, for defining our automorphic forms. This
will take care of the level-raising arguments needed later on.  Let $S_r$ be a finite set of places of $F^+$, split in $F$ and disjoint
from $S_\ell$. For $v\in S_r$, let $U_{0,v}\subset G(F^+_v)$ be a compact open subgroup, and let \[ \chi_v:U_{0,v}\to\Or^\times \] be a
morphism with open kernel. We will use the notation $U_{r}=\prod_{v\in S_r}U_{0,v}$ and $\chi=\prod_{v\in S_r}\chi_v$.

Fix the data $\{\ba,U_r,\chi\}$. Let $M_{\ba,\chi}=M_{\ba}\otimes_\Or(\bigotimes_{v\in S_r}\Or(\chi_v))$. Let $U\subset
G(\A_{F^+}^\infty)$ be a compact open subgroup such that its projection to the $v$-th coordinate is contained in $U_{0,v}$ for each $v\in
S_r$. Let $A$ be an $\Or$-algebra.  Suppose either that the projection of $U$ to $G(F^+_\ell)$ is contained in $G(\Or_{F^+_\ell})$, or that $A$ is
a $K$-algebra. Then define $S_{\ba,\chi}(U,A)$ to be the space of functions \[ f:G(F^+)\backslash G(\A_{F^+}^\infty)\to
M_{\ba,\chi}\otimes_\Or A \] such that \[ f(gu)=u_{\ell,S_r}^{-1}f(g)\hspace{4mm}\forall\hspace{1mm}g\in G(\A_{F^+}^\infty),u\in U,\] where
$u_{\ell,S_r}$ denotes the product of the projections to the coordinates of $S_\ell$ and $S_r$. Here, $u_{S_r}$ acts already on $M_{\ba,\chi}$
by $\chi$, and the action of $u_\ell$ is via $\xi_\ba$.

Let $V$ be any compact subgroup of $G(\A_{F^+}^\infty)$ such that its projection to $G(F^+_v)$ is contained in $U_{0,v}$ for each $v\in S_r$,
and let $A$ be an $\Or$-algebra. If either $A$ is a $K$-algebra or the projection of $V$ to $G(F^+_\ell)$ is contained in
$G(\Or_{F^+_\ell})$, denote by $S_{\ba,\chi}(V,A)$ the union of the $S_{\ba,\chi}(U,A)$, where $U$ runs over compact open subgroups
containing $V$ for which their projection to $G(F^+_v)$ is contained in $U_{0,v}$ for each $v\in S_r$, and for which their projection to
$G(F^+_\ell)$ is contained in $G(\Or_{F^+_\ell})$ if $A$ is not a $K$-algebra.  Note that if $V\subset V'$ then $S_{\ba,\chi}(V',A)\subset
S_{\ba,\chi}(V,A)$.

If $U$ is open and we choose a decomposition \[ G(\A_{F^+}^\infty)=\coprod_{j\in J}G(F^+)g_jU,\] then the map $f\mapsto(f(g_j))_{j\in J}$
defines an injection of $A$-modules \begin{equation}\label{suffsmall} S_{\ba,\chi}(U,A)\hookrightarrow\prod_{j\in J}M_{\ba,\chi}\otimes_\Or A.\end{equation} Since $G(F^+)\backslash G(\A_{F^+}^\infty)/U$ is finite and $M_{\ba,\chi}$ is a free $\Or$-module of finite rank, we have that
$S_{\ba,\chi}(U,A)$ is a finitely generated $A$-module.

We say that a compact open subgroup $U\subset G(\A_{F^+}^\infty)$ is sufficiently small if for some finite place $v$ of $F^+$, the projection of $U$ to $G(F^+_v)$ contains only one element of finite order. Note that the map (\ref{suffsmall}) is not always surjective, but it is if, for example,
$U$ is sufficiently small. Thus, in this case, $S_{\ba,\chi}(U,A)$ is a free $A$-module of rank 
\[ (\dim_{K}W_{\ba}).\#\left(G(F^+)\backslash G(\A_{F^+}^\infty)/U\right).\]
 Moreover,
if either $U$ is sufficiently small or $A$ is $\Or$-flat, we have that \[ S_{\ba,\chi}(U,A)=S_{\ba,\chi}(U,\Or)\otimes_\Or A. \]

Let $U$ and $V$ be compact subgroups of $G(\A_{F^+}^\infty)$ such that their projections to $G(F^+_v)$ are contained in $U_{0,v}$ for each
$v\in S_r$. Suppose either $A$ is a $K$-algebra or that the projections of $U$ and $V$ to $G(F^+_\ell)$ are contained in $G(\Or_{F^+_\ell})$.
Also, let $g\in G(\A_{F^+}^{S_r,\infty})\times U_r$; if $A$ is not a $K$-algebra, we suppose that $g_\ell\in G(\Or_{F^+_\ell})$. If
$V\subset gUg^{-1}$, then there is a natural map \[ g:S_{\ba,\chi}(U,A)\longrightarrow S_{\ba,\chi}(V,A) \] defined by \[
(gf)(h)=g_{\ell,S_r}f(hg).\] In particular, if $V$ is a normal subgroup of $U$, then $U$ acts on $S_{\ba,\chi}(V,A)$, and we have that \[
S_{\ba,\chi}(U,A)=S_{\ba,\chi}(V,A)^U.\]

Let $U_1$ and $U_2$ be compact subgroups of $G(\A_{F^+}^\infty)$ such that their projections to $G(F^+_v)$ are contained in $U_{0,v}$ for all
$v\in S_r$. Let $g\in G(\A_{F^+}^{S_r,\infty})\times U_r$. If $A$ is not a $K$-algebra, we suppose that the projections of $U_1$ and $U_2$
to $G(F^+_\ell)$ are contained in $G(\Or_{F^+_\ell})$, and that $g_\ell\in G(\Or_{F^+_\ell})$. Suppose also that $\#U_1gU_2/U_2<\infty$ (this
will be automatic if $U_1$ and $U_2$ are open). Then we can define an $A$-linear map \[ [U_1gU_2]:S_{\ba,\chi}(U_2,A)\longrightarrow
S_{\ba,\chi}(U_1,A) \] by \[ ([U_1gU_2]f)(h)=\sum_i(g_i)_{\ell,S_r}f(hg_i),\] if $U_1gU_2=\coprod_ig_iU_2$.

\begin{lemma}\label{lemacoinv} Let $U\subset G(\A_{F^+}^{\infty,S_r})\times\prod_{v\in S_r}U_{0,v}$ be a sufficiently small compact open subgroup and let $V\subset U$ be a normal open subgroup. Let $A$ be an $\Or$-algebra. Suppose that either $A$ is a $K$-algebra or the projection of $U$ to $G(F^+_\ell)$ is contained in $G(\Or_{F^+_\ell})$. Then $S_{\ba,\chi}(V,A)$ is a finite free $A[U/V]$-module. Moreover, let $I_{U/V}\subset A[U/V]$ be the augmentation ideal and let $S_{\ba,\chi}(V,A)_{U/V}=S_{\ba,\chi}(V,A)/I_{U/V}S_{\ba,\chi}(V,A)$ be the module of coinvariants. Define
\[ \tr_{U/V}:S_{\ba,\chi}(V,A)_{U/V}\to S_{\ba,\chi}(U,A)=S_{\ba,\chi}(V,A)^U \]
as $\tr_{U/V}(f)=\sum_{u\in U/V}uf$. Then $\tr_{U/V}$ is an isomorphism.
\begin{proof} This is the analog of Lemma 3.3.1 of \cite{cht}, and can be proved in the same way.
\end{proof}
\end{lemma}

Choose an isomorphism $\iota:\overline K\overset\simeq\longrightarrow\C$. The choice of $\widetilde I_\ell$ gives a bijection
\begin{equation}\label{iotaplus} \iota_*^+:(\Z^{n,+})^{\Hom(F,K)}_c\overset{\sim}{\longrightarrow}(\Z^{n,+})^{\Hom(F^+,\R)}, \end{equation}
where $(\Z^{n,+})^{\Hom(F,K)}_c$ denotes the set of elements $\ba\in(\Z^{n,+})^{\Hom(F,K)}$ such that
\[ a_{\tau c,i}=-a_{\tau,n+1-i} \]
for every $\tau\in\Hom(F,K)$ and every $i=1,\dots,n$. The map is given by $(\iota_*^+\ba)_{\tau}=a_{\widetilde{\iota^{-1}\tau}}$. We have an isomorphism $\theta:W_\ba\otimes_{K,\iota}\C\rightarrow
W_{\iota_*^+\ba,}$.  Then the map \[ S_{\ba,\emptyset}(\{1\},\C)\longrightarrow S_{(\iota_*^+\ba)^\vee} \] given by \[
f\mapsto(g\mapsto\theta(g_\ell f(g))) \] is an isomorphism of $\C[G(\A_{F^+}^\infty)]$-modules, where, $(\iota_*^+\ba)^\vee_{\tau,i}=-(\iota_*^+\ba)_{\tau,n+1-i}$. Its inverse is given by \[ \phi\mapsto(g\mapsto
g_\ell^{-1}\theta^{-1}(\phi(g))).\] It follows that $S_{\ba,\emptyset}(\{1\},\C)$ is a semi-simple admissible module. Hence, $S_{\ba,\emptyset}(\{1\},\overline K)$ is also
semi-simple admissible, and this easily implies that $S_{\ba,\chi}(U_r,\overline K)$ is a semi-simple admissible
$G(\A_{F^+}^{\infty,S_r})$-module. If $\pi\subset S_{\ba,\emptyset}(\{1\},\overline K)$ is an irreducible $G(\A_{F^+}^{\infty,S_r})\times U_r$-constituent such that the subspace
on which $U_r$ acts by $\chi^{-1}$ is non-zero, then this subspace is an irreducible constituent of $S_{\ba,\chi}(U_r,\overline
K)$, and every irreducible constituent of it is obtained in this way.
 
\subsection{Base change and descent}\label{condnd} Keep the notation as above. We will assume from now on the following hypotheses.
\begin{itemize}
 \item $F/F^+$ is unramified at all finite places.
 \item $G_v$ is quasi-split for every finite place $v$.
\end{itemize}
It is not a very serious restriction for the applications we have in mind, because we will always be able to base change to this situation. First, note that given $F/F^+$, if $n$ is odd there always exists a totally definite unitary group $G$ in $n$ variables with $G_v$ quasi-split for every finite $v$. If $n$ is even, such a $G$ exists if and only if $[F^+:\Q]n/2$ is also even. This follows from the general classification of unitary groups over number fields in terms of the local Hasse invariants.

Let $G_n^*=\Res_{F/F^+}(\GL_n)$. Let $v$ be a
finite place of $F^+$, so that $G_v$ is an unramified group. In particular, it contains hyperspecial maximal compact subgroups. Let $\sigma_{v}$ be any irreducible admissible representation of $G(F_{v}^{+})$. If $v$ is split in $F$, or if $v$ is inert and $\sigma_v$ is spherical, there exists an irreducible admissible representation $\BC_v(\sigma_v)$ of $G_n^*(F^+_v)$, called the {\em local base change} of $\sigma_v$, with the following properties. Suppose that $v$ is inert and $\sigma_v$ is a spherical representation of $G(F^+_v)$; then $\BC_v(\sigma_v)$ is an unramified representation of $G_n^*(F^+_v)$,
whose Satake parameters are explicitly determined in terms of those of $\sigma_v$; the formula is given in \cite{minguez}, where we take the {\em standard base change} defined there. If $v$ splits in $F$ as $ww^c$, the local base change in this case is $\BC_v(\sigma_v)=\sigma_v\circ
i_w^{-1}\otimes(\sigma_v\circ i_{w^c}^{-1})^\vee$ as a representation of $G_n^*(F^+_v)=\GL_n(F_w)\times\GL_n(F_{w^c})$. In this way, if we see $\BC_v(\sigma_v)$ as a representation of $G(F^+_v)\times G(F^+_v)$ via the isomorphism $i_w\times i_{w^c}:G(F^+_v)\times G(F^+_v)\overset{\sim}{\to}\GL_n(F_w)\times\GL_n(F_{w^c})$, then $\BC_v(\sigma_v)=\sigma_v\otimes\sigma_v^\vee$. The base change for ramified finite places is being treated in the work of M{\oe}glin, but for our applications it is enough to assume that $F/F^+$ is unramified at finite places. 

In the global case, if $\sigma$ is an automorphic representation of $G(\A_{F^+})$, we say that an automorphic representation $\Pi$ of
$G_n^*(\A_{F^+})=\GL_n(\A_F)$ is a (strong) base change of $\sigma$ if $\Pi_v$ is the local base change of $\sigma_v$ for every finite $v$, except those inert $v$ where $\sigma_v$ is not spherical, and if the infinitesimal character of $\Pi_\infty$ is the base change of that of $\sigma_\infty$. In particular, since $G(F^+_\infty)$ is compact, $\Pi$ is cohomological.

The following theorem is one of the main results of \cite{lab}, and a key ingredient in this paper. We use the notation $\boxplus$ for the isobaric sum of discrete automorphic representations, as in \cite{clozel}.

\begin{thm}[Labesse]\label{labesse} Let $\sigma$ be an automorphic representation of $G(\A_{F^+})$. Then there exists a partition
\[ n=n_1+\cdots+n_r \]
and discrete, conjugate self dual automorphic representations $\Pi_i$ of $\GL_{n_i}(\A_F)$, for $i=1,\dots,r$, such that
\[ \Pi_1\boxplus\cdots\boxplus\Pi_r \]
is a base change of $\sigma$.

Conversely, let $\Pi$ be a conjugate self dual, cuspidal, cohomological automorphic representation of $\GL_n(\A_F)$. Then $\Pi$ is the base change of an automorphic representation $\sigma$ of $G(\A_{F^+})$. Moreover, if such a $\sigma$ satisfies that $\sigma_v$ is spherical for every inert place $v$ of $F^+$, then $\sigma$ appears with multiplicity one in the cuspidal spectrum of $G$.
\begin{proof} The first part is Corollaire 5.3 of \cite{lab} and the second is Th\'eor\`{e}me 5.4.
\end{proof}
\end{thm}

\begin{rems*}\label{remlabesse}
\begin{enumerate}[(1)]
 \item In \cite{lab} there are two hypothesis to Corollaire 5.3, namely, the property called (*) by Labesse and that $\sigma_\infty$ is a discrete series, which are automatically satisfied in our case because the group is totally definite. 
 \item Since $\Pi_1\boxplus\dots\boxplus\Pi_r$ is a base change of $\sigma$, it is a cohomological representation of $\GL_n(\A_F)$. However, this doesn't imply that each $\Pi_i$ is cohomological, although it will be if $n-n_i$ is even.
 \item The partition $n=n_1+\cdots n_r$ and the representations $\Pi_i$ are uniquely determined by multiplicity one for $\GL_n$, because the $\Pi_i$ are discrete.
\end{enumerate}
\end{rems*}

\subsection{Galois representations of unitary type via unitary groups} Keep the notation and assumptions as in the last sections. 

\begin{thm} Let $\pi$ be as above. Let $\pi=\otimes_{v\not\in S_r}\pi_v$ be an irreducible constituent of the space $S_{\ba,\chi}(U_r,\overline K)$. Then there exists a unique continuous semisimple representation
\[ r_\ell(\pi):\Gal(\overline F/F)\to\GL_n(\overline K) \]
satisfying the following properties.

\begin{enumerate}[(i)]
\item If $v\not\in S_\ell\cup S_r$ is a place of $F^+$ which splits as $v=ww^c$ in $F$, then
\[ r_\ell(\pi)|_{\Gamma_w}^{\ssp}\simeq\left(r_\ell(\pi_v\circ i_w^{-1})\right)^{\ssp}.\]
\item $r_\ell(\pi)^c\cong r_\ell(\pi)^\vee(1-n)$.
\item If $v$ is an inert place such that $\pi_v$ is spherical then $r_\ell(\pi)$ is unramified at $v$.
\item If $w|\ell$ then $r_\ell(\pi)$ is de Rham at $w$, and if moreover $\pi_{w|_{F^+}}$ is unramified, then $r_\ell(\pi)$ is crystalline at $w$.
\item For every $\tau\in\Hom(F,K)$ giving rise to an place $w|\ell$ of $F$, the Hodge-Tate weights of $r|_{\Gamma_w}$ with respect to $\tau$ are given by
\[ \HT_{\tau}(r|_{\Gamma_w})=\{j-n-a_{\tau,j}\}_{j=1,\dots,n}.\]
In particular, $r$ is Hodge-Tate regular.
\end{enumerate}

\begin{proof}
For the uniqueness, note that the set of places $w$ of $F$ which are split over a place $v$ of $F^+$ which is not in $S_\ell\cup S_r$ has Dirichlet density $1$, and hence, if two continuous semisimple representations $\Gal(\overline F/F)\to\GL_n(\overline\Q_\ell)$ satisfy property (i), they are isomorphic.

Take an isomorphism $\iota:\overline K\overset{\sim}{\longrightarrow}\C$. By the above argument, the representation we will construct will not depend on it. By means of $\iota$ and the choice of $\widetilde I_\ell$, we obtain a (necessarily cuspidal) automorphic representation $\sigma=\otimes_v\sigma_v$ of $G(\A_{F^+})$, such that $\sigma_v=\iota\pi_v$ for $v\not\in S_r$ finite and $\sigma_\infty$ is the representation of $G(F^+_\infty)$ given by the weight $(\iota_*^+\ba)^\vee\in(\Z^{n,+})^{\Hom(F^+,\R)}$. By Theorem \ref{labesse}, there is a partition $n=n_1+\dots+n_r$ and discrete automorphic representations $\Pi_i$ of $\GL_{n_i}(\A_F)$ such that
\[ \Pi=\Pi_1\boxplus\cdots\boxplus\Pi_r \]
is a strong base change of $\sigma$. Moreover, $\Pi$ is cohomological of weight $\iota_*\ba$, where $(\iota_*\ba)_\tau=\ba_{\iota^{-1}\tau}$ for $\tau\in\Hom(F,\C)$. For each $i=1,\dots,r$, let $S_i\supset S_\ell$ be any finite set of finite primes of $F^+$, unramified in $F$. For each $i=1,\dots,r$, let $\psi_i:\A_F^\times/F^\times\to\C^\times$ be a character such that
\begin{itemize}
\item $\psi_i^{-1}=\psi_i^c$;
\item $\psi_i$ is unramified above $S_i$, and
\item for every $\tau\in\Hom(F,\C)$ giving rise to an infinite place $w$, we have
\[ \psi_{i,w}(z)=\left(\tau z/|\tau z|\right)^{\delta_{i,\tau}},\]
where $|z|^2=z\overline z$ and $\delta_{i,\tau}=0$ if $n-n_i$ is even, and $\delta_{i,\tau}=\pm1$ otherwise.
\end{itemize}
Thus, if $n-n_i$ is even, we may just choose $\psi_i=1$. The proof of the existence of such a character follows from a similar argument used in the proof of \cite[Lemma VII.2.8]{ht}. With these choices, it follows that $\Pi_i\psi_i$ is cohomological. Also, by the classification of M{\oe}glin and Waldspurger (\cite{mw}), there is a factorization $n_i=a_ib_i$, and a cuspidal automorphic representation $\rho_i$ of $\GL_{a_i}(\A_F)$ such that
\[ \Pi_i\psi_i=\rho_i\boxplus\rho_i\vabs\boxplus\dots\boxplus\rho_i\vabs^{b_i-1}.\]
Moreover, $\rho_i\vabs^{\frac{b_i-1}{2}}$ is cuspidal and conjugate self dual. Let $\chi_i:\A_F^\times/F^\times\to\C^\times$ be a character such that
\begin{itemize}
\item $\chi_i^{-1}=\chi_i^c$;
\item $\chi_i$ is unramified above $S_i$, and
\item for every $\tau\in\Hom(F,\C)$ giving rise to an infinite place $w$, we have
\[ \chi_{i,w}(z)=\left(\tau z/|\tau z|\right)^{\mu_{i,\tau}},\]
where $\mu_{i,\tau}=0$ if $a_i$ is odd or $b_i$ is odd, and $\mu_{i,\tau}=\pm1$ otherwise.
\end{itemize}

Then $\rho_i\vabs^{\frac{b_i-1}{2}}\chi_i$ is cuspidal, cohomological and conjugate self dual. Note that $\chi_i^{-1}\vabs^{(a_i-1)(b_i-1)/2}$ and $\psi_i^{-1}\vabs^{\frac{n_i-n}{2}}$ are algebraic characters. Let
\[
\begin{array}{rl} r_\ell(\pi) = & \bigoplus_{i=1}^r\left(r_\ell\left(\rho_i\chi_i\vabs^{\frac{b_i-1}{2}}\right)\otimes\epsilon^{a_i-n_i}\otimes r_\ell\left(\chi_i^{-1}\vabs^{(a_i-1)(b_i-1)/2}\right)\right. \\ & \left.\otimes\left(1\oplus\epsilon\oplus\cdots\oplus\epsilon^{b_i-1}\right)\otimes r_\ell\left(\psi_i^{-1}\vabs^{\frac{n_i-n}{2}}\right)\right),\end{array} \]
where $r_\ell=r_{\ell,\iota}$ and $\epsilon$ is the $\ell$-adic cyclotomic character. This is a continuous semisimple representation which satisfies all the required properties. We use the freedom to vary the sets $S_i$ to achieve property (iii).
\end{proof}
\end{thm}

\begin{rem}\label{autounit}
In the proof of the above theorem, if $r=1$ and $\Pi$ is already cuspidal, then $r_{\ell}(\pi)\cong r_{\ell,\iota}(\Pi)$. As a consequence, suppose that $\iota:\overline\Q_\ell\overset{\sim}{\longrightarrow}\C$ is an isomorphism and $\Pi$ is a conjugate self dual, cohomological, cuspidal automorphic representation of $\GL_n(\A_F)$ of weight $\iota_*\ba$. Then, by Theorem \ref{labesse}, we can find an irreducible constituent $\pi\subset S_{\ba,\emptyset}(\{1\},\overline K)$ such that $r_{\ell,\iota}(\Pi)\cong r_{\ell}(\pi)$.
\end{rem}

\begin{rem}\label{resirr} If $r_{\ell}(\pi)$ is irreducible, then the base change of $\pi$ is already cuspidal. Indeed, from the construction made in the proof and Remark \ref{remlabesse}, (2), we see that $r_{\ell}(\pi)$ is a direct sum of $r$ representations $r_i$ of dimension $n_i$. If $r_{\ell}(\pi)$ is irreducible, we must have $r=1$. Similarly, the discrete base change $\Pi$ must be cuspidal, because otherwise there would be a factorization $n=ab$ with $a,b>1$ and $r_{\ell}(\pi)$ would be a direct sum of $b$ representations of dimension $a$. This proves our claim.
\end{rem}

\section{An $R^{\red}=T$ theorem for Hecke algebras of unitary groups}

\subsection{Hecke algebras} Keep the notation and assumptions as in the last section.  For each place $w$ of $F$, split
above a place $v$ of $F^+$, let $\Iw(w)\subset G(\Or_{F^+_v})$ be the inverse image under $i_w$ of the group of matrices in
$\GL_n(\Or_{F_w})$ which reduce modulo $w$ to an upper triangular matrix. Let $\Iw_1(w)$ be the kernel of the natural surjection
$\Iw(w)\to(k_w^\times)^n$, where $k_w$ is the residue field of $F_w$. Similarly, let $U_0(w)$ (resp. $U_1(w)$) be the inverse image
under $i_w$ of the group of matrices in $\GL_n(\Or_{F_w})$ whose reduction modulo $w$ has last row $(0,\dots,0,*)$ (resp.
$(0,\dots,0,1)$). Then $U_1(w)$ is a normal subgroup of $U_0(w)$, and the quotient $U_0(w)/U_1(w)$ is naturally isomorphic to
$k_w^\times$.

Let $Q$ be a finite (possibly empty) set of places of $F^+$ split in $F$, disjoint from $S_\ell$ and $S_r$, and let $T\supset S_r\cup S_\ell\cup Q$ be a finite set of places of $F^+$ split in $F$. Let $\widetilde T$ denote a set of primes of $F$ above $T$ such that $\widetilde T\coprod\widetilde
T^c$ is the set of all primes of $F$ above $T$. For $v\in T$, we denote by $\widetilde v$ the corresponding element of $\widetilde T$, and for $S\subset T$, we denote by $\widetilde S$ the set of places of $F$ consisting of the $\widetilde v$ for $v\in T$. Let
\[ U=\prod_vU_v\subset G(\A_{F^+}^\infty) \]
be a sufficiently small compact open subgroup such that:

\begin{itemize}
\item if $v\not\in T$ splits in $F$ then $U_v=G(\Or_{F^+_v})$;
\item if $v\in S_r$ then $U_v=\Iw(\widetilde v)$;
\item if $v\in Q$ then $U_v=U_1(\widetilde v)$;
\item if $v\in S_\ell$ then $U_v\subset G(\Or_{F^+_v})$.
\end{itemize}
We write $U_r=\prod_{v\in S_r}U_v$. For $v\in S_r$, let $\chi_v$ be an $\Or$-valued character of $\Iw(\widetilde v)$, trivial on $\Iw_1(\widetilde v)$. Since $\Iw(\widetilde v)/\Iw_1(\widetilde v)\simeq(k_{\widetilde v}^\times)^n$, $\chi_v$ is of the form \[ g\mapsto\prod_{i=1}^n\chi_{v,i}(g_{ii}),\] where $\chi_{v,i}:k_{\widetilde v}^\times\to\Or^\times$.

Let $w$ be a place of $F$, split over a place $v$ of $F^+$ which is not in $T$. We translate the Hecke operators $T_{F_w}^{(j)}$ for $j=1,\dots,n$ on $\GL_n(\Or_{F_w})$ to $G$ via the isomorphism $i_w$. More precisely, let $g_w^{(j)}$ denote the element of $G(\A_{F^+}^\infty)$ whose $v$-coordinate is \[ i_w^{-1}\left(
	\begin{array}{cc}\overline\omega_w1_j&0\\ 0&1_{n-j}\end{array} \right),\] and with all other coordinates equal to $1$. Then we
		define $T_w^{(j)}$ to be the operator $[Ug_w^{(j)}U]$ of $S_{\ba,\chi}(U,A)$. We will denote by $\T^T_{\ba,\chi}(U)$
		the $\Or$-subalgebra of $\End_\Or(S_{\ba,\chi}(U,\Or))$ generated by the operators $T_w^{(j)}$ for $j=1,\dots,n$ and
		$(T_w^{(n)})^{-1}$, where $w$ runs over places of $F$ which are split over a place of $F^+$ not in $T$. The algebra $\T^T_{\ba,\chi}(U)$ is reduced, and finite free as an $\Or$-module (see \cite{cht}). Since $\Or$ is a domain, this also implies that $\T^T_{\ba,\chi}(U)$ is a semi-local ring. If $v\in Q$, we can also translate the Hecke operators $V_{\alpha,1}$ of Section 1, for $\alpha\in F_{\widetilde v}^\times$ with non-negative valuation, in exactly the same manner to operators in $S_{\ba,\chi}(U,A)$, and similarly for $V_{\alpha,0}$ if $U_v=U_0(\widetilde v)$.

Write
\begin{equation}\label{decompS}
S_{\ba,\chi}(U,\overline K)=\oplus_\pi\pi^{U}, 
\end{equation}
where $\pi$ runs over the irreducible constituents of $S_{\ba,\chi}(U_r,\overline K)$ for which $\pi^{U}\neq 0$. The Hecke algebra $\T^{T}_{\ba,\chi}(U)$ acts on each $\pi^{U}$ by a scalar, say, by
\[ \lambda_\pi:\T^{T}_{\ba,\chi}(U)\longrightarrow\overline K. \]
Then, $\ker(\lambda_\pi)$ is a minimal prime ideal of $\T^{T}_{\ba,\chi}(U)$, and every minimal prime is of this form. If $\m\subset\T^{T}_{\ba,\chi}(U)$ is a maximal ideal, then
\[ S_{\ba,\chi}(U,\overline K)_\m\neq 0,\]
and localizing at $\m$ kills all the representations $\pi$ such that $\ker(\lambda_\pi)\not\subset\m$. Note also that $\T^{T}_{\ba,\chi}(U)/\m$ is a finite extension of $k$. For $w$ a place of $F$, split over a place $v\not\in T$, we will denote by $\mathbf{T}_w$ the $n$-tuple $(T_w^{(1)},\dots,T_w^{(n)})$ of elements of $\T^{T}_{\ba,\chi}(U)$. We denote by $\overline{\mathbf{T}}_w$ its reduction modulo $\m$. We use the notation of section 2.4.1 of \cite{cht} regarding torsion crystalline representations and Fontaine-Laffaille modules.

\begin{prop} Suppose that $\m$ is a maximal ideal of $\T^{T}_{\ba,\chi}(U)$ with residue field $k$. Then there is a unique continuous semisimple representation \[
	\overline r_{\m}:\Gal(\overline F/F)\to\GL_n(k) \] with the following properties. The first two already
	characterize $\overline r_{\m}$ uniquely.

	\begin{enumerate}[(i)]
\item $\overline r_{\m}$ is unramified at all but finitely many places.
\item If a place $v\not\in T$ splits as $ww^c$ in $F$ then $\overline r_{\m}$ is unramified at $w$ and $\overline r_{\m}(\Frob_w)$ has characteristic
			polynomial $P_{q_w,\overline{\mathbf{T}}_w}(X)$.
\item $\overline r_{\m}^c\cong\overline r_{\m}^\vee(1-n)$.
\item If a place $v$ of $F^+$ is inert in $F$ and if $U_v$ is a hyperspecial maximal compact subgroup of $G(F^+_v)$, then $\overline r_{\m}$ is unramified above $v$.
\item If $w\in\widetilde S_\ell$ is unramified over $\ell$, $U_{w|_{F^+}}=G(\Or_{F^+_w})$ and for every $\tau\in\widetilde I_\ell$ above $w$ we have that
\[ \ell-1-n\geq a_{\tau,1}\geq\dots\geq a_{\tau,n}\geq 0,\]
then
\[ \overline r_{\m}|_{\Gamma_w}=\mathbf{G}_w(\overline M_{\m,w})\]
for some object $\overline M_{\m,w}$ of $\mathscr{MF}_{k,w}$. Moreover, for every $\tau\in\widetilde I_\ell$ over $w$, we have
\[ \dim_{k}(\gr^{-i}\overline M_{\m,w})\otimes_{\Or_{F_w}\otimes_{\Z_\ell}\Or,\tau\otimes 1}\Or=1\]
if $i=j-n-a_{\tau,j}$ for some $j=1,\dots,n$, and $0$ otherwise.
\end{enumerate}
\begin{proof} Choose a minimal prime ideal $\mathfrak{p}\subset{\m}$ and an irreducible constituent $\pi$ of 
\[ S_{\ba,\chi}(U_r,\overline K) \]
such that $\pi^{U}\neq 0$ and $\T^{T}_{\ba,\chi}(U)$ acts on $\pi^{U}$ via $\T^{T}_{\ba,\chi}(U)/\mathfrak{p}$. Choose an invariant lattice for $r_\ell(\pi)$ and define then $\overline r_{\m}$ to be the semi-simplification of the reduction of $r_\ell(\pi)$. This satisfies all of the statements of the proposition, except for the fact that a priori it takes values on the algebraic closure of $k$. Since all the characteristic polynomials of the elements on the image of $\overline r_\m$ have coefficients in $k$, we may assume (because $k$ is finite) that, after conjugation, $\overline r_{\m}$ actually takes values in $k$.
\end{proof}
\end{prop}

We say that a maximal ideal $\mathfrak{m}\subset\T^{T}_{a,\chi}(U)$ is {\em Eisenstein} if $\overline r_\m$ is absolutely reducible. We define (see Chapter 2 of \cite{cht}) $\mathscr{G}_n$ as the group scheme over $\Z$ given by the semi-direct product of $\GL_n\times\GL_1$ by the group $\{1,\j\}$ acting on $\GL_n\times\GL_1$ by 
\[ \j(g,\mu)\j^{-1}=(\mu^tg^{-1},\mu).\]
There is a homomorphism $\nu:\mathscr{G}_n\to\GL_1$ which sends $(g,\mu)$ to $\mu$ and $\j$ to $-1$.

\begin{prop}\label{galnoneis} Let ${\m}$ be a non-Eisenstein maximal ideal of $\T^{T}_{a,\chi}(U)$, with residue field equal to $k$. Then $\overline r_{\m}$
	has an extension to a continuous morphism \[ \overline r_{\m}:\Gal(\overline F/F^+)\to\G_n(k).\] Pick such an extension. Then
	there is a unique continuous lifting \[ r_{\m}:\Gal(\overline F/F^+)\to\G_n(\T^{T}_{a,\chi}(U)_\m) \] of $\overline r_{\m}$ with the
	following properties. The first two of these already characterize the lifting $r_{\m}$ uniquely.

\begin{enumerate}[(i)] \item $r_{\m}$ is unramified at almost all places.
                  \item If a place $v\not\in T$ of $F^+$ splits as $ww^c$ in $F$, then $r_{\m}$ is unramified at $w$ and $r_{\m}(\Frob_w)$ has characteristic polynomial $P_{q_w,\mathbf{T}_w}(X)$.
		\item $\nu\circ r_{\m}=\epsilon^{1-n}\delta_{F/F^+}^{\mu_\m}$, where $\delta_{F/F^+}$ is the non-trivial character of $\Gal(F/F^+)$
		and $\mu_\m\in\Z/2\Z$.  \item If $v$ is an inert place of $F^+$ such that $U_v$ is a hyperspecial maximal compact subgroup of $G(F^+_v)$ then $r_{\m}$ is unramified at $v$.
		\item Suppose that $w\in\widetilde S_\ell$ is unramified over $\ell$, that $U_{w|_{F^+}}=G(\Or_{F^+_w})$, and that for every $\tau\in\widetilde I_\ell$ above $w$ we have that
		\[ \ell-1-n\geq a_{\tau,1}\geq\dots\geq a_{\tau,n}\geq 0.\]
		Then for each open ideal $I\subset\T^T_{\ba,\chi}(U)_\m$, 
		\[ \left(r_{\m}\otimes_{\T^T_{\ba,\chi}(U)_\m}\T^T_{\ba,\chi}(U)_\m/I\right)|_{\Gamma_w}=\mathbf{G}_w(M_{\m,I,w}) \]
		for some object $M_{\m,I,w}$ of $\mathscr{MF}_{\Or,w}$.
		
		 \item If $v\in S_r$ and $\sigma\in I_{F_{\widetilde v}}$ then $r_{\m}(\sigma)$ has characteristic polynomial \[
	 	\prod_{j=1}^n(X-\chi_{v,j}^{-1}(\Art_{F_{\widetilde v}}^{-1}\sigma)).\]

 \item Suppose that $v\in Q$. Let $\phi_{\widetilde v}$ be a lift of $\Frob_{\widetilde v}$ to $\Gal(\overline F_{\widetilde v}/F_{\widetilde v})$. Suppose that $\alpha\in k$ is a simple root of the characteristic polynomial of $\overline r_{\m}(\phi_{\widetilde v})$. Then there exists a unique root $\widetilde\alpha\in\T^{T}_{\ba,\chi}(U)_{\m}$ of the characteristic polynomial of $r_{\m}(\phi_{\widetilde v})$ which lifts $\alpha$.

Let $\overline\omega_{\widetilde v}$ be the uniformizer of $F_{\widetilde v}$ corresponding to $\phi_{\widetilde v}$ via $\Art_{F_{\widetilde v}}$. Suppose that $Y\subset S_{\ba,\chi}(U,K)_{\m}$ is a $\T^{T}_{\ba,\chi}(U)[V_{\varpi_{\widetilde v},1}]$-invariant subspace such that $V_{\varpi_{\widetilde v},1}-\widetilde\alpha$ is topologically nilpotent on $Y$, and let $\T^{T}(Y)$ denote the image of $\T^{T}_{\ba,\chi}(U)$ in $\End_{K}(Y)$. Then for each $\beta\in F_{\widetilde v}^\times$ with non-negative valuation, $V_{\beta,1}$ (in $\End_{K}(Y)$) lies in $\T^{T}(Y)$, and $\beta\mapsto V(\beta)$ extends to a continuous character $V:F_{\widetilde v}^\times\to\T^{T}(Y)^{\times}$. Further, $(X-V_{\varpi_{\widetilde v},1})$ divides the characteristic polynomial of $r_{\m}(\phi_{\widetilde v})$ over $\T^{T}(Y)$.

Finally, if $q_v\equiv 1\mod\ell$ then
\[ r_{\m}|_{\Gamma_{\widetilde v}}\cong s\oplus(V\circ\Art_{F_{\widetilde v}}^{-1}), \]
where $s$ is unramified.

\end{enumerate}
\begin{proof} This is the analogue of Proposition 3.4.4 of \cite{cht}, and can be proved exactly in the same way.
\end{proof}
\end{prop}

\begin{coro}\label{coroTw} Let $Q'$ denote a finite set of places of $F^+$, split in $F$ and disjoint from $T$. Let $\m$ be a non-Eisenstein maximal ideal of $\T^T_{\ba,\chi}(U)$ with residue field $k$, and let $U_1(Q')=\prod_{v\not\in Q'}U_v\times\prod_{v\in Q'}U_1(\widetilde v)$. Denote by $\varphi:\T^{T\cup Q'}_{\ba,\chi}(U')\to\T^T_{\ba,\chi}(U)$ the natural map, and let $\m'=\varphi^{-1}(\m)$, so that $\m'$ is also non-Eisenstein with residue field $k$. Then the localized map $\varphi:\T^{T\cup Q'}_{\ba,\chi}(U_{1}(Q'))_{\m'}\to\T^T_{\ba,\chi}(U)_{\m}$ is surjective.
\begin{proof} It suffices to see that $T_w^{(j)}/1$ is in the image of $\varphi$ for $j=1,\dots,n$ and $w$ a place of $F$ over $Q'$, which follows easily because $r_\m=\varphi\circ r_{\m'}$, and so
\[ T_w^{(j)}=\varphi\left(q_w^{j(1-j)/2}\tr\left(\bigwedge^jr_{\m'}\right)(\phi_w)\right), \]
where $\phi_{w}$ is any lift of Frobenius at $w$.
\end{proof}
\end{coro}

\subsection{The main theorem} In this section we will use the Taylor-Wiles method in the version improved by Diamond, Fujiwara, Kisin and Taylor. We will recapitulate the running assumptions made until now, and add a few more. Thus, let $F^+$ be a totally real field and $F/F^+$ a totally imaginary quadratic extension. Fix a positive integer $n$ and an odd prime $\ell>n$. Let $K/\Q_\ell$ be a finite extension, let $\overline K$ be an algebraic closure of $K$, and suppose that $K$ is big enough to contain the image of every embedding $F\hookrightarrow\overline K$. Let $\Or$ be the ring of integers of $K$, and $k$ its residue field. Let $S_\ell$ denote the set of places of $F^+$ above $\ell$. Let $\widetilde S_\ell$ denote a set of places of $F$ above $\ell$ such that $\widetilde S_\ell\coprod\widetilde S_\ell^c$ are all the places above $\ell$. We let $\widetilde I_\ell$ denote the set of embeddings $F\hookrightarrow K$ which give rise to a place in $\widetilde S_\ell$. We will suppose that the following conditions are satisfied.

\begin{itemize}
 \item $F/F^+$ is unramified at all finite places;
 \item $\ell$ is unramified in $F^+$;
 \item every place of $S_\ell$ is split in $F$;
\end{itemize}

Let $G$ be a totally definite unitary group in $n$ variables, attached to the extension $F/F^+$ such that $G_v$ is quasi-split for every finite place $v$ (cf. Section \ref{condnd} for conditions on $n$ and $[F^+:\Q]$ to ensure that such a group exists). Choose a lattice in $F^+$ giving a model for $G$ over $\Or_{F^+}$, and fix a basis of the lattice, so that for each split $v=ww^c$, there are two isomorphisms
\[ i_w:G_v\longrightarrow\GL_{n/F_w} \]
and
\[ i_{w^c}:G_v\longrightarrow\GL_{n/F_{w^c}} \]
taking $G(\Or_{F^+_v})$ to $\GL_n(\Or_{F_w})$ and $\GL_n(\Or_{F_{w^c}})$ respectively.

Let $S_a$ denote a finite, non-empty set of primes of $F^+$, disjoint from $S_\ell$, such that if $v\in S_a$ then
\begin{itemize}
 \item $v$ splits in $F$, and
 \item if $v$ lies above a rational prime $p$ then $v$ is unramified over $p$ and $[F(\zeta_p):F]>n$.
\end{itemize}

Let $S_r$ denote a finite set of places of $F^+$, disjoint from $S_a\cup S_\ell$, such that if $v\in S_r$ then
\begin{itemize}
 \item $v$ splits in $F$, and
 \item $q_v\equiv 1\mod\ell$.
\end{itemize}

We will write $T=S_\ell\cup S_a\cup S_r$, and $\widetilde T\supset\widetilde S_\ell$ for a set of places of $F$ above those of $T$ such that $\widetilde T\coprod\widetilde T^c$ is the set of all places of $F$ above $T$. For $S\subset T$, we will write $\widetilde S$ to denote the set of $\widetilde v$ for $v\in S$. We will fix a compact open subgroup
\[ U=\prod_vU_v \]
of $G(\A_{F^+}^\infty)$, such that
\begin{itemize}
 \item if $v$ is not split in $F$ then $U_v$ is a hyperspecial maximal compact subgroup of $G(F^+_v)$;
 \item if $v\not\in S_a\cup S_r$ splits in $F$ then $U_v=G(\Or_{F^+_v})$;
 \item if $v\in S_r$ then $U_v=\Iw(\widetilde v)$, and
 \item if $v\in S_a$ then $U_v=i_{\widetilde v}^{-1}\ker(\GL_n(\Or_{F_{\widetilde v}})\to\GL_n(k_{\widetilde v}))$.
\end{itemize}
Then, $U$ is sufficiently small ($U_v$ has only one element of finite order if $v\in S_a$) and its projection to $G(F^+_\ell)$ is contained in $G(\Or_{F^+_\ell})$. Write 
\[ U_r=\prod_{v\in S_r}U_v.\]
For any finite set $Q$ of places of $F^+$, split in $F$ and disjoint from $T$, we will write $T(Q)=T\cup Q$. Also, we will fix a set of places $\widetilde T(Q)\supset\widetilde T$ of $F$ over $T(Q)$ as above, for each $Q$. We will also write
\[ U_0(Q)=\prod_{v\not\in Q}U_v\times\prod_{v\in Q}U_0(\widetilde v) \]
and
\[ U_1(Q)=\prod_{v\not\in Q}U_v\times\prod_{v\in Q}U_1(\widetilde v).\]
Thus, $U_0(Q)$ and $U_1(Q)$ are also sufficiently small compact open subgroups of $G(\A_{F^+}^\infty)$.

Fix an element $\ba\in(\Z^{n,+})^{\Hom(F,K)}$ such that for every $\tau\in\widetilde I_\ell$ we have
\begin{itemize}
 \item $a_{\tau c,i}=-a_{n+1-i}$ and
 \item $\ell-1-n\geq a_{\tau,1}\geq\dots\geq a_{\tau,n}\geq 0$.
\end{itemize}

Let $\m\subset\T^T_{\ba,1}(U)$ be a non-Eisenstein maximal ideal with residue field equal to $k$. Write $\T=\T^T_{\ba,1}(U)_\m$. Consider the representation
\[ \overline r_\m:\Gal(\overline F/F^+)\to\G_n(k) \]
and its lifting
\[ r_\m:\Gal(\overline F/F^+)\to\G_n(\T) \]
given by Proposition \ref{galnoneis}. For $v\in T$, denote by $\overline r_{\m,v}$ the restriction of $\overline r_\m$ to a decomposition group $\Gamma_{\widetilde v}$ at $\widetilde v$. We will {\em assume} that $\overline r_\m$ has the following properties.

\begin{itemize}
 \item $\overline r_\m(\Gal(\overline F/F^+(\zeta_\ell)))$ is big (see Definition 2.5.1 of \cite{cht}, where the same notion is also defined for subgroups of $\GL_n(k)$);
 \item if $v\in S_r$ then $\overline r_{\m,v}$ is the trivial representation of $\Gamma_{\widetilde v}$, and
 \item if $v\in S_a$ then $\overline r_\m$ is unramified at $v$ and
 \[ H^0(\Gamma_{\widetilde v},(\ad\overline r_\m)(1))=0.\]
\end{itemize}

We will use the Galois deformation theory developed in Section 2 of \cite{cht}, to where we refer the reader for the definitions and results. Consider the global deformation problem
\[ \CS=(F/F^+,T,\widetilde T,\Or,\overline r_\m,\epsilon^{1-n}\delta_{F/F^+}^{\mu_\m},\{\mathscr{D}_v\}_{v\in T}),\]
where the local deformation problems $\mathscr{D}_v$ are as follows. For $v\in T$, we denote by
\[ r_v^{\univ}:\Gamma_{\widetilde v}\to\GL_n(R_v^{\loc}) \]
the universal lifting ring of $\overline r_{\m,v}$, and by $\mathscr{I}_v\subset R_v^{\loc}$ the ideal corresponding to $\mathscr{D}_v$.
\begin{itemize}
 \item For $v\in S_a$, $\mathscr{D}_v$ consists of all lifts of $\overline r_{\m,v}$, and thus $\mathscr{I}_v=0$.
 \item For $v\in S_\ell$, $\mathscr{D}_v$ consist of all lifts whose artinian quotients all arise from torsion Fontaine-Laffaille modules, as in Section 2.4.1 of \cite{cht}.
 \item For $v\in S_r$, $\mathscr{D}_v$ corresponds to the ideal $\mathscr{I}_v^{(1,1,\dots,1)}$ of $R_v^{\loc}$, as in Section 3 of \cite{taylorII}. Thus, $\mathscr{D}_v$ consists of all the liftings $r:\Gamma_{\widetilde v}\to\GL_n(A)$ such that for every $\sigma$ in the inertia subgroup $I_{\widetilde v}$, the characteristic polynomial of $r(\sigma)$ is
 \[ \prod_{i=1}^n(X-1).\]
\end{itemize}
Let
\[ r_{\CS}^{\univ}:\Gal(\overline F/F^+)\to\G_n(R_{\CS}^{\univ}) \]
denote the universal deformation of $\overline r_\m$ of type $\CS$. By Proposition \ref{galnoneis}, $r_\m$ gives a lifting of $\overline r_\m$ which is of type $\mathscr{S}$; this gives rise to a surjection
\[ R_{\CS}^{\univ}\longrightarrow\T.\]
Let $H=S_{\ba,1}(U,\Or)_{\m}$. This is a $\T$-module, and under the above map, a $R_{\CS}^{\univ}$-module. Our main result is the following.

\begin{thm}\label{mainthm} Keep the notation and assumptions of the start of this section. Then
 \[ (R_{\CS}^{\univ})^{\red}\simeq\T.\]
Moreover, $\mu_\m\equiv n\mod 2$.
\begin{proof}
 The proof is essentially the same as Taylor's (\cite{taylorII}), except that here there are no primes $S(B)_1$ and $S(B)_2$, in his notation. One has just to note that his argument is still valid in our simpler case. The idea is to use Kisin's version (\cite{kisin}) of the Taylor-Wiles method in the following way, in order to avoid dealing with non-minimal deformations separately. There are essentially two moduli problems to consider at places in $S_r$. One of them consists in considering all the characters $\chi_v$ to be trivial. This is the case in which we are ultimately interested, but the local deformation rings are not so well behaved (for example, they are not irreducible). We call this the {\em degenerate case}. On the other hand, we can also consider the characters $\chi_v$ in such a way that $\chi_{v,i}\neq\chi_{v,j}$ for all $v\in S_r$ and all $i\neq j$. This is the {\em non-degenerate case}, and we can always consider such a set of characters by our assumption that $\ell>n$. Note that both problems are equal modulo $\ell$. The Taylor-Wiles-Kisin method doesn't work with the first moduli problem, but it works fine in the non-degenerate case. Taylor's idea is to apply all the steps of the method simultaneosly for the degenerate and non-degenerate cases, and obtain the final conclusion of the theorem by means of comparing both processes modulo $\lambda$, and using the fact that in the degenerate case, even if the local deformation ring is not irreducible, every prime ideal which is minimal over $\lambda$ contains a unique minimal prime, and this suffices to proof what we want. We will reproduce most of the argument in the following pages. What we will prove in the end is that $H$ is a nearly faithful $R_{\CS}^{\univ}$-module, which by definition means that the ideal $\Ann_{R_{\CS}^{\univ}}(H)$ is nilpotent. Since $\T$ is reduced, this proves the main statement of the theorem.

We will be working with several deformation problems at a time. Consider a set $Q$ of finite set of places of $F^+$, disjoint from $T$, such that if $v\in Q$, then
\begin{itemize}
 \item $v$ splits as $ww^c$ in $F$,
 \item $q_v\equiv 1\mod\ell$, and
 \item $\overline r_{\m,v}=\overline\psi_{v}\oplus\overline s_{v}$, with $\dim\overline\psi_{v}=1$ and such that $\overline s_{v}$ does not contain $\overline\psi_{v}$ as a sub-quotient.
\end{itemize}

Let $T(Q)$ and $\widetilde T(Q)$ be as in the start of the section. Also, let $\{\chi_v:\Iw(\widetilde v)/\Iw_1(\widetilde v)\to\Or^\times\}_{v\in S_r}$ be a set of characters of order dividing $\ell$. To facilitate the notation, we will write $\chi_v=(\chi_{v,1},\dots,\chi_{v,n})$ and $\chi=\{\chi_v\}_{v\in S_r}$. Consider the deformation problem given by
\[ \mathscr{S}_{\chi,Q}=(F/F^+,T(Q),\widetilde T(Q),\Or,\overline r_\m,\epsilon^{1-n}\delta_{F/F^+}^{\mu_\m},\{\mathscr{D}'_v\}_{v\in T(Q)}),
\]
where:
\begin{itemize}
 \item for $v\in S_a\cup S_\ell$, $\mathscr{D}'_v=\mathscr{D}_v$;
 \item for $v\in S_r$, $\mathscr{D}'_v$ consists of all the liftings $r:\Gamma_{\widetilde v}\to\GL_n(A)$ such that the characteristic polynomial of $r(\sigma)$ for $\sigma\in I_{\widetilde v}$ is
 \[ \prod_{i=1}^n(X-\chi_{v,i}^{-1}(\Art_{F_{\widetilde v}}^{-1}\sigma)) \]
 (see Section 3 of \cite{taylorII}).
 \item for $v\in Q$, $\mathscr{D}'_v$ consists of all Taylor-Wiles liftings of $\overline r_{\m,v}$, as in Section 2.4.6 of \cite{cht}. More precisely, $\mathscr{D}'_v$ consists of all the liftings $r:\Gamma_{\widetilde v}\to\GL_n(A)$ which are conjugate to one of the form $\psi_v\oplus s_v$ with $\psi_v$ a lift of $\overline\psi_{v}$ and $s_v$ an unramified lift of $\overline s_{v}$.
\end{itemize}
Denote by $\mathscr{I}_v^{\chi_v}$ the corresponding ideal of $R_v^{\loc}$ for every $v\in T(Q)$. Let
\[ r_{\CS_{\chi,Q}}^{\univ}:\Gal(\overline F/F^+)\to\G_n(R_{\CS_{\chi,Q}}^{\univ}) \]
denote the universal deformation of $\overline r$ of type $\CS_{\chi,Q}$, and let
\[ r_{\CS_{\chi,Q}}^{\square_T}:\Gal(\overline F/F^+)\to\G_n(R_{\CS_{\chi,Q}}^{\square_T}) \]
denote the universal $T$-framed deformation of $\overline r$ of type $\CS_{\chi,Q}$ (see \cite[2.2.7]{cht} for the definition of $T$-framed deformations; note that it depends on $\widetilde T$). Thus, by definition of the deformation problems, we have that $R_{\CS_{1,\emptyset}}^{\univ}=R_{\CS}^{\univ}$. As we claimed above, both problems are equal modulo $\ell$. We have isomorphisms
\begin{equation}\label{runivQ}
 R_{\CS_{\chi,Q}}^{\univ}/\lambda\cong R_{\CS_{1,Q}}^{\univ}/\lambda
\end{equation}
and
\begin{equation}\label{rsquareQ}
 R_{\CS_{\chi,Q}}^{\square_T}/\lambda\cong R_{\CS_{1,Q}}^{\square_T}/\lambda,
\end{equation}
compatible with the natural commutative diagrams
\[\xymatrix{
 R_{\CS_{\chi,Q}}^{\univ}\ar@{->>}[r]\ar[d] & R_{\CS_{\chi,\emptyset}}^{\univ} \ar[d] \\
 R_{\CS_{\chi,Q}}^{\square_T}\ar@{->>}[r] & R_{\CS_{\chi,\emptyset}}^{\square_T} }
\]
and
\[\xymatrix{
 R_{\CS_{1,Q}}^{\univ}\ar@{->>}[r]\ar[d] & R_{\CS_{1,\emptyset}}^{\univ} \ar[d] \\
 R_{\CS_{1,Q}}^{\square_T}\ar@{->>}[r] & R_{\CS_{1,\emptyset}}^{\square_T} }
\]
Also, let
\[ R_{\chi,T}^{\loc}=\widehat\bigotimes_{v\in T}R_v^{\loc}/\mathscr{I}_v^{\chi_v}.
\]
Then
\begin{equation}\label{rlocchi}
 R_{\chi,T}^{\loc}/\lambda\cong R_{1,T}^{\loc}/\lambda.
\end{equation}
To any $T$-framed deformation of type $\mathscr{S}_{\chi,Q}$ and any $v\in T$ we can associate a lifting of $\overline r_{\m,v}$ of type $\mathscr{D}_v$, and hence there are natural maps
\[ R_{\chi,T}^{\loc}\longrightarrow R_{\CS_{\chi,Q}}^{\square_T} \]
which modulo $\lambda$ are compatible with the identifications (\ref{rlocchi}) and (\ref{rsquareQ}).

Let $\mathscr{T}=\Or[[X_{v,i,j}]]_{v\in T;i,j=1,\dots,n}$. Then a choice of a lifting $r_{\CS_{\chi,Q}}^{\univ}$ of $\overline r_\m$ over $R_{\CS_{\chi,Q}}^{\univ}$ representing the universal deformation of type $\CS_{\chi,Q}$ gives rise to an isomorphism of $R_{\CS,Q}^{\univ}$-algebras
\begin{equation}\label{tensor} 
R_{\CS_{\chi,Q}}^{\square_T}\simeq R_{\CS_{\chi,Q}}^{\univ}\hat\otimes_\Or\mathscr{T}, 
\end{equation}
so that
\[ (r_{\CS_{\chi,Q}}^{\univ};\{1_n+(X_{v,i,j})\}_{v\in T}) \]
represents the universal $T$-framed deformation of type $\CS_{\chi,Q}$ (see Proposition 2.2.9 of \cite{cht}). Moreover, we can choose the liftings $r_{\CS_{\chi,Q}}^{\univ}$ so that
\[ r_{\CS_{\chi,Q}}^{\univ}\otimes_\Or k=r_{\CS_{1,Q}}^{\univ}\otimes_\Or k \]
under the natural identifications (\ref{runivQ}). Then the isomorphisms (\ref{tensor}) for $\chi$ and $1$ are compatible with the identifications (\ref{rsquareQ}) and (\ref{runivQ}). 

For $v\in Q$, let $\psi_v$ denote the lifting of $\overline\psi_{\widetilde v}$ to $(R_{\CS_{\chi,Q}}^{\univ})^\times$ given by the lifting $r_{\CS_{\chi,Q}}^{\univ}$. Also, write $\Delta_Q$ for the maximal $\ell$-power order quotient of $\prod_{v\in Q}k_{\widetilde v}^\times$, and let $\mathfrak{a}_Q$ denote the ideal of $\mathscr{T}[\Delta_Q]$ generated by the augmentation ideal of $\Or[\Delta_Q]$ and by the $X_{v,i,j}$ for $v\in T$ and $i,j=1,\dots,n$. Since the primes of $Q$ are different from $\ell$ and $\overline\psi_{\widetilde v}$ is unramified, $\psi_v$ is tamely ramified, and then
\[ \prod_{v\in Q}(\psi_v\circ\Art_{F_{\widetilde v}}):\Delta_Q\longrightarrow(R_{\CS_{\chi,Q}}^{\univ})^\times \]
makes $R_{\CS_{\chi,Q}}^{\univ}$ an $\Or[\Delta_Q]$-algebra. This algebra structure is compatible with the identifications (\ref{runivQ}), because we chose the liftings $r_{\CS_{\chi,Q}}^{\univ}$ and $r_{\CS_{1,Q}}^{\univ}$ compatibly. Via the isomorphisms (\ref{tensor}), $R_{\CS_{\chi,Q}}^{\square_T}$ are $\mathscr{T}[\Delta_Q]$-algebras, which are compatible modulo $\lambda$ for the different choices of $\chi$. Finally, we have an isomorphism
\begin{equation}\label{aug}
 R_{\CS_{\chi,Q}}^{\square_T}/\mathfrak{a}_Q\simeq R_{\CS_{\chi,\emptyset}}^{\univ},
\end{equation}
compatible with the identifications (\ref{rsquareQ}) and (\ref{runivQ}), the last one with $Q=\emptyset$.

Note that since
\[ S_{\ba,1}(U,k)=S_{\ba,\chi}(U,k)\]
we can find a maximal ideal $\m_{\chi,\emptyset}\subset\T^T_{\ba,\chi}(U)$ with residue field $k$ such that for a prime $w$ of $F$ split over a prime $v\not\in T$ of $F^+$, the Hecke operators $T_w^{(j)}$ have the same image in $\T^T_{\ba,\chi}(U)/\m_{\chi,\emptyset}=k$ as in $\T^T_{\ba,1}(U)/\m=k$. It follows that $\overline r_{\m_{\chi,\emptyset}}\cong\overline r_\m$, and in particular $\m_{\chi,\emptyset}$ is non-Eisenstein. We define $\m_{\chi,Q}\subset\T^{T(Q)}_{\ba,\chi}(U_1(Q))$ as the preimage of $\m_{\chi,\emptyset}$ under the natural map
\[ \T^{T(Q)}_{\ba,\chi}(U_1(Q))\twoheadrightarrow\T^{T(Q)}_{\ba,\chi}(U_0(Q))\twoheadrightarrow\T^{T(Q)}_{\ba,\chi}(U)\hookrightarrow\T^T_{\ba,\chi}(U).\]
Then $\T^{T(Q)}_{\ba,\chi}(U_1(Q))/\m_{\chi,Q}=k$, and if a prime $w$ of $F$ splits over a prime $v\not\in T(Q)$ of $F^+$, then the Hecke operators $T_w^{(j)}$ have the same image in $\T^{T(Q)}_{\ba,\chi}(U_1(Q))/\m_{\chi,Q}=k$ as in $\T^{T}_{\ba,1}(U)/\m=k$. Hence, $\overline r_{\m_{\chi,Q}}\cong\overline r_\m$ and $\m_{\chi,Q}$ is non-Eisenstein. Let
\[ r_{\m_{\chi,Q}}:\Gal(\overline F/F^+)\to\G_n(\T^{T(Q)}_{\ba,\chi}(U_{1}(Q))_{\m_{\chi,Q}}) \]
be the continuous representation attached to $\m_{\chi,Q}$ as in Proposition \ref{galnoneis}. Write $\T_\chi=\T^T_{\ba,\chi}(U)_{\m_{\chi,\emptyset}}$ and $H_\chi=S_{\ba,\chi}(U,\Or)_{\m_{\chi,\emptyset}}$. We have the following natural surjections
\begin{equation}\label{surjT}
\T^{T(Q)}_{\ba,\chi}(U_{1}(Q))_{\m_{\chi,Q}}\twoheadrightarrow\T^{T(Q)}_{\ba,\chi}(U_{0}(Q))_{\m_{\chi,Q}}\twoheadrightarrow\T^{T(Q)}_{\ba,\chi}(U)_{\m_{\chi,Q}}=\T_\chi.
\end{equation}
The last equality follows easily from Corollary \ref{coroTw}.

For each $v\in Q$, choose $\phi_{\widetilde v}\in \Gamma_{\widetilde v}$ a lift of $\Frob_{\widetilde v}$, and let $\overline\omega_{\widetilde v}\in F_{\widetilde v}^\times$ be the uniformizer corresponding to $\phi_{\widetilde v}$ via $\Art_{F_{\widetilde v}}$.  Let
\[ P_{\widetilde v}\in\T^{T(Q)}_{\ba,\chi}(U_{1}(Q))_{\m_{\chi,Q}}[X] \]
denote the characteristic polynomial of $r_{\m_{\chi,Q}}(\phi_{\widetilde v})$. Since $\overline\psi_{v}(\phi_{\widetilde v})$ is a simple root of the characteristic polynomial of $\overline r_\m(\phi_{\widetilde v})$, by Hensel's lemma, there exists a unique root $A_{\widetilde v}\in\T^{T(Q)}_{\ba,\chi}(U_{1}(Q))_{\m_{\chi,Q}}$ of $P_{\widetilde v}$ lifting $\overline\psi_v(\phi_{\widetilde v})$. Thus, there is a factorisation
\[ P_{\widetilde v}(X)=(X-A_{\widetilde v})Q_{\widetilde v}(X) \]
over $\T^{T(Q)}_{\ba,\chi}(U_{1}(Q))_{\m_{\chi,Q}}$, where $Q_{\widetilde v}(A_{\widetilde v})\in\T^{T(Q)}_{\ba,\chi}(U_{1}(Q))_{\m_{\chi,Q}}^\times$. By part (i) of Proposition 1.7 and Lemma 1.9, $P_{\widetilde v}(V_{\varpi_{\widetilde v},1})=0$ on $S_{\ba,\chi}(U_{1}(Q),\Or)_{\m_{\chi,Q}}$. For $i=0,1$, let 
\[
H_{i,\chi,Q}=\left(\prod_{v\in Q}Q_{\widetilde v}(V_{\varpi_{\widetilde v},i})\right)S_{\ba,\chi}(U_{i}(Q),\Or)_{\m_{\chi,Q}}\subset S_{\ba,\chi}(U_{i}(Q),\Or)_{\m_{\chi,Q}},\]
and let $\T_{i,\chi,Q}$ denote the image of $\T_{\ba,\chi}^{T(Q)}(U_{1}(Q))_{\m_{\chi,Q}}$ in $\End_{\Or}(H_{i,\chi,Q})$. We see that $H_{1,\chi,Q}$ is a direct summand of $S_{\ba,\chi}(U_1(Q),\Or)$ as a $\T^{T(Q)}_{\ba,\chi}(U_1(Q))$-module. Also, we have an isomorphism
\[ \left(\prod_{v\in Q}Q_{\widetilde v}(V_{\overline\omega_{\widetilde v},0})\right):H_\chi\cong H_{0,\chi,Q}.\]
This can be proved using Proposition \ref{eigenV} and Lemmas \ref{lemma315} and \ref{321}, as in \cite[3.2.2]{cht}.

For all $v\in Q$, $V_{\varpi_{\widetilde v},1}=A_{\widetilde v}$ on $H_{1,\chi,Q}$. By part (vii) of Proposition \ref{galnoneis}, for each $v\in Q$ there is a character with open kernel
\[ V_{v}:F_{\widetilde v}^\times\longrightarrow\T_{1,\chi,Q}^\times \]
such that
\begin{itemize}
 \item if $\alpha\in\Or_{F_{\widetilde v}}$ is non-zero, then $V_{\alpha,1}=V_{v}(\alpha)$ on $H_{1,\chi,Q}$ and
 \item $(r_{\m_{\chi,Q}}\otimes\T_{1,\chi,Q})|_{\Gamma_{\widetilde v}}\cong s_v\oplus(V_{v}\circ\Art_{F_{\widetilde v}}^{-1})$, where $s_v$ is unramified.
\end{itemize}
It is clear that $V_v\circ\Art_{F_{\widetilde v}}^{-1}$ is a lifting of $\overline\psi_v$ and $s_v$ is a lifting of $\overline s_v$. It follows by (v) and (vi) of the same proposition that $r_{\m_{\chi,Q}}\otimes\T_{1,\chi,Q}$ gives rise to a deformation of $\overline r_\m$ of type $\CS_{\chi,Q}$, and thus to a surjection
\[ R_{\CS_{\chi,Q}}^{\univ}\twoheadrightarrow\T_{1,\chi,Q},\]
such that the composition
\[ \prod_{v\in Q}\Or_{F_{\widetilde v}}^\times\twoheadrightarrow\Delta_Q\to(R_{\CS_{\chi,Q}}^{\univ})^\times\to\T_{1,\chi,Q}^\times \]
coincides with $\prod_{v\in Q}V_{v}$. We then have that $H_{1,\chi,Q}$ is an $R_{\CS_{\chi,Q}}^{\univ}$-module, and we set
\[ H_{1,\chi,Q}^{\square_T}=H_{1,\chi,Q}\otimes_{R_{\CS_{\chi,Q}}^{\univ}}R_{\CS_{\chi,Q}}^{\square_T}=H_{1,\chi,Q}\otimes_\Or\mathscr{T}.\]

Since $\ker(\prod_{v\in Q}k_{\widetilde v}^\times\to\Delta_Q)$ acts trivially on $H_{1,\chi,Q}$ and $H_{1,\chi,Q}$ is a $\T^{T(Q)}_{\ba,\chi}(U_{1}(Q))$-direct summand of $S_{\ba,\chi}(U_{1}(Q),\Or)$, Lemma \ref{lemacoinv} implies that $H_{1,\chi,Q}$ is a finite free $\Or[\Delta_Q]$-module, and that
\[ (H_{1,\chi,Q})_{\Delta_Q}\cong H_{0,\chi,Q}\cong H_\chi.\]

Since $U$ is sufficiently small, we get isomorphisms
\[ S_{\ba,\chi}(U,\Or)\otimes_{\Or}k\cong S_{\ba,\chi}(U,k)=S_{\ba,1}(U,k)\cong S_{\ba,1}(U,\Or)\otimes_{\Or}k \]
and
\[ S_{\ba,\chi}(U_{1}(Q),\Or)\otimes_{\Or}k\cong S_{\ba,\chi}(U_{1}(Q),k)=S_{\ba,1}(U_{1}(Q),k)\cong S_{\ba,1}(U_{1}(Q),\Or)\otimes_{\Or}k. \]
Thus we get identifications
\[ H_{\chi}/\lambda\cong H_{1}/\lambda, \]
\[ H_{1,\chi,Q}/\lambda\cong H_{1,1,Q}/\lambda \]
and
\[ H_{1,\chi,Q}^{\square_T}/\lambda\cong H_{1,1,Q}^{\square_T}/\lambda,\]
compatible with all the pertinent identifications modulo $\lambda$ made before. 

Let
\[ \varepsilon_\infty=(1-(-1)^{\mu_\m-n})/2 \]
and
\[ q_0=[F^+:\Q]n(n-1)/2+[F^+:\Q]n\varepsilon_\infty.\]
By Proposition 2.5.9 of \cite{cht}, there is an integer $q\geq q_0$, such that for every natural number $N$, we can find a set of primes $Q_N$ (and a set of corresponding $\overline\psi_v$ and $\overline s_v$ for $\overline r_\m$) such that
\begin{itemize}
 \item $\#Q_N=q$;
 \item for $v\in Q_N$, $q_v\equiv 1(\Mod\ell^N)$ and
 \item $R_{\CS_{1,Q_N}}^{\square_T}$ can be topologically generated over $R_{1,T}^{\loc}$ by $q'=q-q_0$ elements.
\end{itemize}
Define
\[ R_{\chi,\infty}^{\square_T}=R_{\chi,T}^{\loc}[[Y_1,\dots,Y_{q'}]].\]
Then there is a surjection
\[ R_{1,\infty}^{\square_T}\twoheadrightarrow R_{\CS_{1,Q_N}}^{\square_T} \]
extending the natural map $R_{1,T}^{\loc}\to R_{\CS_{1,Q_N}}^{\square_T}$. Reducing modulo $\lambda$ and lifting the obtained surjection, via the identifications
\[ R_{\chi,\infty}^{\square_T}/\lambda\simeq R_{1,\infty}^{\square_T}/\lambda, \]
we obtain a surjection
\[ R_{\chi,\infty}^{\square_T}\twoheadrightarrow R_{\CS_{\chi,Q_N}}^{\square_T} \]
extending the natural map $R_{\chi,T}^{\loc}\to R_{\CS_{\chi,Q_N}}^{\square_T}$.

For $v\in S_a$, $R_v^{\loc}/\mathscr{I}_v^{\chi_v}$ is a power series ring over $\Or$ in $n^2$ variables (see Lemma 2.4.9 of \cite{cht}), and for $v\in S_\ell$ it is a power series ring over $\Or$ in $n^2+[F_{\widetilde v}:\Q_\ell]n(n-1)/2$ variables (see Corollary 2.4.3 of {\em loc. cit.}).

Suppose that $\chi_{v,i}\neq\chi_{v,j}$ for every $v\in S_r$ and every $i,j=1,\dots,n$ with $i\neq j$. Then, by Proposition 3.1 of \cite{taylorII}, for every $v\in S_r$, $R_v^{\loc}/\mathscr{I}_v^{\chi_v}$ is irreducible of dimension $n^2+1$ and its generic point has characteristic zero. It follows that $(R_v^{\loc}/\mathscr{I}_v^{\chi_v})^{\red}$ is geometrically integral (in the sense that $(R_v^{\loc}/\mathscr{I}_v^{\chi_v})^{\red}\otimes_{\Or}\Or'$ is an integral domain for every finite extension $K'/K$, where $\Or'$ is the ring of integers of $K'$) and flat over $\Or$. Moreover, by part 3. of Lemma 3.3 of \cite{cy2},
\[ (R_{\chi,\infty}^{\square_T})^{\red}\simeq\left(\left(\widehat\bigotimes_{v\in S_r}(R_v^{\loc}/\mathscr{I}_v^{\chi_v})^{\red}\right)\widehat\bigotimes\left(\widehat\bigotimes_{v\in S_a\cup S_\ell}R_v^{\loc}/\mathscr{I}_v\right)\right)[[Y_1,\dots,Y_{q'}]],\]
and the same part of that lemma implies that $(R_{\chi,\infty}^{\square_T})^{\red}$ is geometrically integral. We conclude that in the non-degenerate case, $R_{\chi,\infty}^{\square_T}$ is irreducible, and, by part 2., its Krull dimension is
\[ 1+q+n^2\#T-[F^+:\Q]n\varepsilon_\infty. \]

Suppose now that we are in the degenerate case, that is, $\chi_v=1$ for every $v\in S_r$. Then (see Proposition 3.1 of \cite{taylorII}) for every such $v$, $R_v^{\loc}/\mathscr{I}_v^{\chi_v}$ is pure of dimension $n^2+1$, its generic points have characteristic zero, and every prime of $R_v^{\loc}/\mathscr{I}_v^{\chi_v}$ which is minimal over $\lambda(R_v^{\loc}/\mathscr{I}_v^{\chi_v})$ contains a unique minimal prime. After eventually replacing $K$ by a finite extension $K'$ (which we are allowed to do since the main theorem for one $K$ implies the same theorem for every $K'$), $R_v^{\loc}/\mathscr{I}_v^{\chi_v}$ satisfies that for every prime ideal $\mathfrak{p}$ which is minimal (resp. every prime ideal $\mathfrak{q}$ which is minimal over $\lambda(R_v^{\loc}/\mathscr{I}_v^{\chi_v})$), the quotient $(R_v^{\loc}/\mathscr{I}_v^{\chi_v})/\mathfrak{p}$ (resp. $(R_v^{\loc}/\mathscr{I}_v^{\chi_v})/\mathfrak{q}$) is geometrically integral. It follows then by parts 2., 5. and 7. of Lemma 3.3 of \cite{cy2} that every prime ideal of $R_{1,\infty}^{\square_T}$ which is minimal over $\lambda R_{1,\infty}^{\square_T}$ contains a unique minimal prime, the generic points of $R_{1,\infty}^{\square_T}$ have characteristic zero and $R_{1,\infty}^{\square_T}$ is pure.

Let $\Delta_\infty=\Z_\ell^q$, $S_\infty=\mathscr{T}[[\Delta_\infty]]$ and $\mathfrak{a}=\ker(S_\infty\twoheadrightarrow\Or)$, where the map sends $\Delta_\infty$ to $1$ and the variables $X_{v,i,j}$ to $0$. Thus, $S_\infty$ is isomorphic to a power series ring over $\Or$ in $q+n^2\#T$ variables. For every $N$, choose a surjection
\[ \Delta_\infty\twoheadrightarrow\Delta_{Q_N}. \]
We have an induced map on completed group algebras
\[ \Or[[\Delta_\infty]]\twoheadrightarrow\Or[\Delta_{Q_N}].\]
and thus a map
\begin{equation}\label{mapatch} S_\infty\twoheadrightarrow\mathscr{T}[\Delta_{Q_N}]\to R_{\CS_{\chi,Q_N}}^{\square_T} \end{equation}
which makes $R_{\CS_{\chi,Q_N}}^{\square_T}$ an algebra over $S_\infty$. The map $S_\infty\twoheadrightarrow\mathscr{T}[\Delta_{Q_N}]$ sends the ideal $\mathfrak{a}$ to $\mathfrak{a}_{Q_N}$. Let $\mathfrak{c}_N=\ker(S_\infty\twoheadrightarrow\mathscr{T}[\Delta_{Q_N}])$. Note that every open ideal of $S_\infty$ contains $\mathfrak{c}_N$ for some $N$. The following properties hold.

\begin{itemize}
 \item $H_{1,\chi,Q_N}^{\square_T}$ is finite free over $S_\infty/\mathfrak{c}_N$.
 \item $R_{\CS_{\chi,Q_N}}^{\square_T}/\mathfrak{a}\simeq R_{\CS_{\chi,\emptyset}}^{\univ}$.
 \item $H_{1,\chi,Q_N}^{\square_T}/\mathfrak{a}\simeq H_\chi$.
\end{itemize}

In what follows, we will use that we can patch the $R_{\CS_{\chi,Q_N}}^{\square_T}$ to obtain in the limit a copy of $R_{\chi,\infty}^{\square_T}$, and simultaneously patch the $H_{1,\chi,Q_N}$ to form a module over $R_{\chi,\infty}^{\square_T}$, finite free over $S_\infty$. The patching construction is carried on in exactly the same way as in \cite{taylorII}. The outcome of this process is a family of $R_{\chi,\infty}^{\square_T}\widehat\otimes_\Or S_\infty$-modules $H_{1,\chi,\infty}^{\square_T}$ with the following properties.
 \begin{enumerate}[(1)]
  \item They are finite free over $S_\infty$, and the $S_\infty$-action factors through $R_{\chi,\infty}^{\square_T}$, in such a way that the obtained maps $S_\infty\to R_{\chi,\infty}^{\square_T}\twoheadrightarrow R_{\CS_{\chi,Q_N}}^{\square_T}$ are the maps defined in (\ref{mapatch}) for every $N$; in particular, there is a surjection 
  \[ R_{\chi,\infty}^{\square_T}/\mathfrak{a}\twoheadrightarrow R_{\CS_{\chi,Q_N}}^{\univ}/\mathfrak{a}=R_{\CS_{\chi,\emptyset}}^{\univ}.\]  
  \item There are isomorphism $H_{1,\chi,\infty}^{\square_T}/\lambda\simeq H_{1,1,\infty}^{\square_T}/\lambda$ of $R_{\chi,\infty}^{\square_T}/\lambda\simeq R_{1,\infty}^{\square_T}/\lambda$-modules.
  \item There are isomorphisms $H_{1,\chi,\infty}^{\square_T}/\mathfrak{a}\simeq H_\chi$ of $R_{\chi,\infty}^{\square_T}/\mathfrak{a}$-modules, where we see $H_\chi$ as a module over $R_{\chi,\infty}^{\square_T}/\mathfrak{a}$ by means of the map in (1). Moreover, these isomorphisms agree modulo $\lambda$ via the identifications of (2).
 \end{enumerate}

Let us place ourselves in the non-degenerate case. That is, let us choose the characters $\chi$ such that $\chi_{v,i}\neq\chi_{v,j}$ for every $v\in S_r$ and every $i\neq j$.  This is possible because $\ell>n$ and $q_v\equiv 1(\Mod\ell)$ for $v\in S_r$. Since the action of $S_\infty$ on $H_{1,\chi,\infty}^{\square_T}$ factors through $R_{\chi,\infty}^{\square_T}$,
\begin{equation}\label{dep1}
\depth_{R_{\chi,\infty}^{\square_T}}(H_{1,\chi,\infty}^{\square_T})\geq\depth_{S_\infty}(H_{1,\chi,\infty}^{\square_T}).
\end{equation}
Also, since $H_{1,\chi,\infty}^{\square_T}$ is finite free over $S_\infty$, which is a Cohen-Macaulay ring, by the Auslander-Buchsbaum formula we have that
\begin{equation}\label{dep2}
\depth_{S_\infty}(H_{1,\chi,\infty}^{\square_T})=\dim S_\infty=1+q+n^2\#T.
\end{equation}
Since the depth of a module is at most its Krull dimension, by equations (\ref{dep1}) and (\ref{dep2}) we obtain that
\begin{equation}\label{depth} \dim\left(R_{\chi,\infty}^{\square_T}/\Ann_{R_{\chi,\infty}^{\square_T}}(H_{1,\chi,\infty}^{\square_T})\right) \geq 1+q+n^2\#T.
\end{equation}
Recall that $R_{\chi,\infty}^{\square_T}$ is irreducible of dimension 
\begin{equation}\label{dim}
1+q+n^2\#T-[F^+:\Q]n\varepsilon_\infty. 
\end{equation}
Then, (\ref{depth}), (\ref{dim}) and Lemma 2.3 of \cite{taylorII} imply that $\varepsilon=0$ (that is, $\mu_\m\equiv n(\Mod 2)$) and that $H_{1,\chi,\infty}^{\square_T}$ is a nearly faithful $R_{\chi,\infty}^{\square_T}$-module. This implies in turn that $H_{1,\chi,\infty}^{\square_T}/\lambda\simeq H_{1,1,\infty}^{\square_T}/\lambda$ is a nearly faithful $R_{\chi,\infty}^{\square_T}/\lambda\simeq R_{1,\infty}^{\square_T}/\lambda$-module (this follows from Nakayama's Lemma, as in Lemma 2.2 of \cite{taylorII}). Since the generic points of $R_{1,\infty}^{\square_T}$ have characteristic zero, $R_{1,\infty}^{\square_T}$ is pure and every prime of $R_{1,\infty}^{\square_T}$ which is minimal over $\lambda R_{1,\infty}^{\square_T}$ contains a unique minimal prime of $R_{1,\infty}^{\square_T}$, the same lemma implies that $H_{1,1,\infty}^{\square_T}$ is a nearly faithful $R_{1,\infty}^{\square_T}$-module. Finally, using the same Lemma again, this implies that $H_{1,1,\infty}^{\square_T}/\mathfrak{a}\simeq H$ is a nearly faithful $R_{1,\infty}^{\square_T}/\mathfrak{a}$-module, and since $R_{1,\infty}^{\square_T}/\mathfrak{a}\twoheadrightarrow R_{\CS}^{\univ}$, $H$ is a nearly faithful $R_{\CS}^{\univ}$-module.
\end{proof}

\end{thm}

\section{Modularity lifting theorems}\label{secmodthm}
In this section we apply the results of the previous sections to prove modularity lifting theorems for $\GL_n$. We deal first with the case of a totally imaginary field $F$.

\begin{thm} Let $F^+$ be a totally real field, and $F$ a totally imaginary quadratic extension of $F^+$. Let $n\geq 1$ be an integer and $\ell>n$ be a prime number, unramified in $F$. Let
 \[ r:\Gal(\overline F/F)\longrightarrow\GL_n(\overline\Q_\ell) \]
be a continuous irreducible representation with the following properties. Let $\overline r$ denote the semisimplification of the reduction of $r$.

\begin{enumerate}[(i)]
 \item $r^c\cong r^\vee(1-n)$.
 \item $r$ is unramified at all but finitely many primes.
 \item For every place $v|\ell$ of $F$, $r|_{\Gamma_v}$ is crystalline.
 \item There is an element $\ba\in(\Z^{n,+})^{\Hom(F,\overline\Q_\ell)}$ such that

     \begin{itemize} 
         \item for all $\tau\in\Hom(F^+,\overline\Q_\ell)$, we have either
		\[ \ell-1-n\geq a_{\tau,1}\geq\dots\geq a_{\tau,n}\geq 0 \]
		or
		\[ \ell-1-n\geq a_{\tau c,1}\geq\dots\geq a_{\tau c,n}\geq 0;\]
	\item for all $\tau\in\Hom(F,\overline\Q_\ell)$ and every $i=1,\dots,n$,
    		 \[ a_{\tau c,i}=-a_{\tau,n+1-i}. \]
	\item for all $\tau\in\Hom(F,\overline\Q_\ell)$ giving rise to a prime $w|\ell$,
	\[ \HT_{\tau}(r|_{\Gamma_w})=\{j-n-a_{\tau,j}\}_{j=1}^n.\]
	In particular, $r$ is Hodge-Tate regular.
	\end{itemize}
 \item $\overline F^{\ker(\ad\overline r)}$ does not contain $F(\zeta_\ell)$.
 \item The group $\overline r(\Gal(\overline F/F(\zeta_\ell)))$ is big.
 \item The representation $\overline r$ is irreducible and there is a conjugate self-dual, cohomological, cuspidal automorphic representation $\Pi$ of $\GL_{n}(\A_{F})$, of weight $\ba$ and unramified above $\ell$, and an isomorphism $\iota:\overline\Q_\ell\overset{\sim}{\longrightarrow}\C$, such that $\overline r\cong\overline r_{\ell,\iota}(\Pi)$.
\end{enumerate}
Then $r$ is automorphic of weight $\ba$ and level prime to $\ell$.

\begin{proof}

Arguing as in \cite[Theorem 5.2]{taylorII}, we may assume that $F$ contains an imaginary quadratic field $E$ with an embedding $\tau_E:E\hookrightarrow\overline\Q_\ell$ such that
\[ \ell-1-n\geq a_{\tau,1}\geq\dots\geq a_{\tau,n}\geq 0\]
for every $\tau:F\hookrightarrow\overline\Q_\ell$ extending $\tau_E$. This will allow us to choose the set $\widetilde S_\ell$ (in the notation of Section \ref{sell}) in such a way that the weights $a_{\tau,i}$ are all within the correct range for $\tau\in\widetilde I_\ell$. Let $\iota:\overline\Q_\ell\overset\simeq\longrightarrow\C$ and let $\Pi$ be a conjugate self dual, cuspidal, cohomological automorphic representation of $\GL_{n}(\A_{F})$ of weight $\iota_*\ba$, with $\Pi_\ell$ unramified, such that $\overline r\cong\overline r_{\ell,\iota}(\Pi)$. Let $S_r$ denote the places of $F$ not dividing $\ell$ at which $r$ or $\Pi$ is ramified. Since $\overline F^{\ker(\ad\overline r)}$ does not contain $F(\zeta_\ell)$, we can choose a prime $v_1$ of $F$ with the following properties.

\begin{itemize}
 \item $v_1\not\in S_r$ and $v_1\nmid\ell$.
 \item $v_1$ is unramified over a rational prime $p$, for which $[F(\zeta_p):F]>n$.
 \item $v_1$ does not split completely in $F(\zeta_\ell)$.
 \item $\ad\overline r(\Frob_{v_1})=1$.
\end{itemize}

Choose a totally real field $L^+/F^+$ with the following properties.

\begin{itemize}
 \item $2|[L^+:\Q]$.
 \item $L^+/F^+$ is Galois and soluble.
 \item $L=L^+E$ is unramified over $L^+$ at every finite place.
 \item $L$ is linearly disjoint from $\overline F^{\ker(\overline r)}(\zeta_\ell)$ over $F$.
 \item $\ell$ is unramified in $L$.
 \item All primes of $L$ above $S_r\cup\{v_1\}$ are split over $L^+$.
 \item $v_1$ and $cv_1$ split completely in $L/F$.
 \item Let $\Pi_L$ denote the base change of $\Pi$ to $L$. If $v$ is a place of $L$ above $S_r$, then 
  \begin{itemize}
    \item $Nv\equiv 1(\Mod\ell)$;
    \item $\overline r(\Gal(\overline L_v/L_v))=1$;
    \item $r|_{I_v}^{\ssp}=1$, and
    \item $\Pi_{L,v}^{\Iw(v)}\neq 0$.
  \end{itemize}
\end{itemize}

Since $[L^+:\Q]$ is even, there exists a unitary group $G$ in $n$ variables attached to $L/L^+$ which is totally definite and such that $G_v$ is quasi-split for every finite place $v$ of $L^+$. Let $S_\ell(L^+)$ denote the set of primes of $L^+$ above $\ell$, $S_r(L^+)$ the set of primes of $L^+$ lying above the restriction to $F^+$ of an element of $S_r$, and $S_a(L^+)$ the set of primes of $L^+$ above $v_1|_{F^+}$. Let $T(L^+)=S_\ell(L^+)\cup S_r(L^+)\cup S_a(L^+)$. It follows from Remarks \ref{autounit} and \ref{resirr} and Theorem \ref{mainthm} that $r|_{\Gal(\overline F/L)}$ is automorphic of weight $\ba_L$ and level prime to $\ell$, where $\ba_L\in(\Z^{n,+})^{\Hom(L,\overline\Q_\ell)}$ is defined as $\ba_{L,\tau}=\ba_{\tau|_F}$. By Lemma 1.4 of \cite{cy2} (note that the hypotheses there must say ``$r^\vee\cong r^c\otimes\chi$'' rather than ``$r^\vee\cong r\otimes\chi$''), this implies that $r$ itself is automorphic of weight $\ba$ and level prime to $\ell$.

\end{proof}
\end{thm}

We can also prove a modularity lifting theorem for totally real fields $F^+$. The proof goes exactly like that of Theorem 5.4 of \cite{taylorII}, using Lemma 1.5 of \cite{cy2} instead of Lemma 4.3.3 of \cite{cht}.

\begin{thm} Let $F^+$ be a totally real field. Let $n\geq 1$ be an integer and $\ell>n$ be a prime number, unramified in $F$. Let
 \[ r:\Gal(\overline F^+/F^+)\longrightarrow\GL_n(\overline\Q_\ell) \]
be a continuous irreducible representation with the following properties. Let $\overline r$ denote the semisimplification of the reduction of $r$.

\begin{enumerate}[(i)]
 \item $r^\vee\cong r({n-1})\otimes\chi$ for some character $\chi:\Gal(\overline F^+/F^+)\to\overline\Q_\ell^\times$ with $\chi(c_v)$ independent of $v|\infty$ (here $c_v$ denotes a complex conjugation at $v$).
 \item $r$ is unramified at all but finitely many primes.
 \item For every place $v|\ell$ of $F$, $r|_{\Gamma_v}$ is crystalline.
 \item There is an element $\ba\in(\Z^{n,+})^{\Hom(F^+,\overline\Q_\ell)}$ such that

     \begin{itemize} 
         \item for all $\tau\in\Hom(F^+,\overline\Q_\ell)$, we have either
		\[ \ell-1-n\geq a_{\tau,1}\geq\dots\geq a_{\tau,n}\geq 0 \]
		or
		\[ \ell-1-n\geq a_{\tau c,1}\geq\dots\geq a_{\tau c,n}\geq 0;\]
	\item for all $\tau\in\Hom(F^+,\overline\Q_\ell)$ and every $i=1,\dots,n$,
    		 \[ a_{\tau c,i}=-a_{\tau,n+1-i}. \]
	\item for all $\tau\in\Hom(F^+,\overline\Q_\ell)$ giving rise to a prime $v|\ell$,
	\[ \HT_{\tau}(r|_{\Gamma_v})=\{j-n-a_{\tau,j}\}_{j=1}^n.\]
	In particular, $r$ is Hodge-Tate regular.
	\end{itemize}
 \item $(\overline F^+)^{\ker(\ad\overline r)}$ does not contain $F^+(\zeta_\ell)$.
 \item The group $\overline r(\Gal(\overline F^+/F^+(\zeta_\ell)))$ is big.
 \item The representation $\overline r$ is irreducible and automorphic of weight $\ba$.
\end{enumerate}
Then $r$ is automorphic of weight $\ba$ and level prime to $\ell$.
\end{thm}

\nocite{ParisBook}

\bibliography{modunitc.bib}{} \bibliographystyle{amsalpha}

\end{document}